\documentclass[preprint,12pt,authoryear,number]{elsarticle} 

\journal{Pure \& Applied Mathematics}

\setcitestyle{square}

\begin{document}

\begin{frontmatter}
\title{Contracting dynamical systems in Banach spaces}

\author[inst1]{Anand Srinivasan\corref{cor1}}
\ead{asrini@alum.mit.edu}
\affiliation[inst1]{organization={Department of Mathematics, MIT},
            addressline={77 Massachusetts Avenue}, 
            city={Cambridge},
            postcode={02139}, 
            state={MA},
            country={USA}
            }

\author[inst2]{Jean-Jacques Slotine}
\ead{jjs@mit.edu}
\affiliation[inst2]{organization={Nonlinear Systems Laboratory, MIT},
            addressline={77 Massachusetts Avenue}, 
            city={Cambridge},
            postcode={02139}, 
            state={MA},
            country={USA}
            }
            
\cortext[cor1]{Corresponding author}

\begin{abstract}
Contraction rates of time-varying maps induced by dynamical systems illuminate a wide range of asymptotic properties with applications in stability analysis and control theory. In finite-dimensional smoothly varying inner-product spaces such as $\R^n$ and $\C^n$ with Riemannian metrics, contraction rates can be estimated by upper-bounding the real numerical range of the vector field's Jacobian. However, vector spaces with norms other than $L^2$ commonly arise in the stability analysis of infinite-dimensional systems such as those arising from partial differential equations and continuum mechanics. To this end, we present a unified approach to contraction analysis in Banach spaces using the theory of weighted semi-inner products. We generalize contraction in a geodesic distance to asymptotic stability of perturbations in smoothly varying semi-inner products, and show that the latter is a dynamical invariant similar to the coordinate-invariance of Lyapunov exponents. We show that contraction in particular weighted spaces verifies asymptotic convergence to subspaces and submanifolds, present applications to limit-cycle analysis and phase-locking phenomena, and pose general conditions for inheritance of contraction properties within coupled systems. We discuss contraction rates in Sobolev spaces for retention of regularity in partial differential equations, and suggest a type of weak solution defined by a vanishing contractive term. Lastly, we present an application to machine learning, using weighted semi-inner products to derive stability conditions for functional gradient descent in a Banach space.
\end{abstract}




\end{frontmatter}

\section{Introduction}
A \textit{contracting} dynamical system in a metric space $(X, d)$ is one for whom any two trajectories $x_1(t), x_2(t)$ satisfy the growth bound:
\begin{align}
    \label{eq:contracting_sys}
    d(x_1(t), x_2(t)) \le e^{\lambda t}d(x_1(0), x_2(0))
\end{align}
for some $\lambda < 0$. The contraction property is both weaker and stronger than traditional notions of Lyapunov and exponential stability in that trajectories uniformly approach one another, but may not approach any fixed point in space. Combined with existence criteria such as fixed points or invariant sets, contraction can be used to recover these standard convergence results.

When $(X, d)$ is a complete metric space, the term \textit{contraction} can be traced back to the Banach fixed-point theorem, where a map $F: X \to X$ which is a \textit{contraction map} with constant $e^\lambda$ sends any point $p$ to an equilibrium $p^*$ with rate $e^{\lambda t}d(p, p_*)$. If $X$ is furthermore a normed vector space, it is known as a Banach space, which is the setting for our analysis of contracting (finite- or infinite-dimensional) systems of differential equations. We denote Banach spaces throughout by $V$ for unambiguity.

In the case of well-posed linear time-invariant (LTI) systems $\dot{u} = Au$ with $A \in \B(V)$ in a Banach space $V$, the contraction map $F$ is given by a one-parameter \textit{strongly continuous} (referred to as $C^0$) semigroup $\Phi(t) : \R^+ \to \B(V)$ which generalizes the matrix exponential $e^{tA}$. In this case, the contraction rate \eqref{eq:contracting_sys} is equivalent to the operator norm bound $\norm{\Phi(t)} \le e^{\lambda t}$ (see \cite{crandall1972nonlinear} and \cite{pazy2012semigroups} for treatment of  linear and nonlinear semigroup theory). As in the finite-dimensional case, spectral theory plays a central role in contraction analysis of linear systems, via the Hille-Yosida theorem which completely characterizes contraction semigroups ($\lambda \le 0$) via the resolvent of the generator $A$.

However, when dynamical systems in Banach spaces are induced by differential equations such as PDEs, moving to an explicit propagator description presents a difficulty in requiring explicit (classical or weak) solutions of the equation(s). Thus, a cornerstone of contraction analysis was the development of sufficient conditions which could be directly applied to a vector field $\dot{u} = f(t, u)$ to verify the exponential convergence of any two solutions $u_1(t), u_2(t)$ in \eqref{eq:contracting_sys}. This argument proceeds by application of the fundamental theorem of calculus and Gr\"onwall's inequality to the growth of perturbations (see \ref{sec:pointfree}). With this approach, we transition from a \textit{point-free} analysis of propagators $\Phi(t, u_0)$ to a \textit{point-wise} analysis of particular trajectories $u(t)$,  in doing so substituting contraction plus \textit{existence} properties (of solutions, equilibria, invariant sets, etc.) for well-posedness requirements. In section \ref{sec:pde}, we give conditions in which the former implies the latter for PDEs, such as uniqueness and regularity properties. 

In the general case of nonlinear nonautonomous systems in $\R^n$ ($\C^n$), \cite{lohmiller1998contraction} showed that the contraction rate in \eqref{eq:contracting_sys} is bounded as
\begin{align}
    \label{eq:l2_contraction}
    \lambda &\le \sup_{v \in V} \frac{\Re \langle v, J(u)v\rangle}{\langle v, v\rangle} = \max\left(\Spec \left(\frac{J(u) + J(u)^\intercal}{2}\right)\right)
\end{align}
where $J(u) = Df(u)$ is the Jacobian of a vector field on $\R^n$ ($\C^n$). One important contribution was the generalization of contraction rates to smoothly time- and space-varying inner-products $\langle \cdot, \cdot\rangle_{t, p}$, thus giving a notion of contraction in a time-varying geodesic distance. Subsequently, \cite{cisneros2020contraction} generalized \eqref{eq:l2_contraction} to infinite-dimensional Hilbert spaces, and \cite{simpson2014contraction} generalized contraction in smoothly-varying inner products on $\R^n$ and $\C^n$ to contraction on Riemannian manifolds.

These contraction results take place in the setting of inner-product spaces. However, while all (separable) Hilbert spaces of fixed dimension are equivalent up to linear isometry, Banach spaces are not and may have norms inducing distinct topologies in the infinite-dimensional setting. This property is critical, for example, in establishing comparison principles and growth estimates of partial differential equations using methods such as Sobolev embedding (see \cite{evans1998partial}). Prior classical and modern works such as \cite{soderlind2006logarithmic, aminzare2014contraction} have extended the quadratic form appearing in the numerator \eqref{eq:l2_contraction} to Banach spaces via numerous other functionals, including \textit{semi-inner products}, a construction on normed spaces originally posed by  \cite{lumer1961semi} and subsequently analyzed by \cite{giles1967classes} for partially extending Hilbert-space notions such as dual pairings and numerical ranges to normed spaces. However, due to the strong norm-dependence of contraction rates (a property we compare and contrast to Lyapunov exponents in Section \ref{sec:metric_dependence}), it is often necessary to find a suitable coordinate transformation which is contracting, which was one of the major contributions of \cite{lohmiller1998contraction} in $\ell^2(\R^n)$. 

To this end, in this work we generalize the theory of contraction in time-varying Riemannian metrics on $\R^n$ to contraction in smoothly-weighted semi-inner products on $V$. We show how this choice of weight can be use to deduce rich properties of dynamical systems in Banach spaces, such as convergence to invariant sets, eventual symmetries, pattern-formation, and the inheritance of contraction properties in coupled systems. 

\section{Background}
\subsection{Contracting dynamical systems in Banach spaces}
\label{sec:evolution_family}
We begin with some background for vector-valued differential equations.
\begin{definition}[Banach space]
A normed vector space $V$ is a \textit{Banach space} if it is a complete metric space in the induced distance
$$
d_V(u, v) = \norm{u - v}
$$
Equivalently, $V$ is a Banach space if any only if every absolutely convergent series is convergent.
\end{definition}
Every Hilbert space is a Banach space, but the norm of the latter does not in general obey the parallelogram law, the geometric property that the norm $\norm{\cdot}_V$ is induced by a Hermitian inner product $\langle\cdot, \cdot\rangle_V$. Thus, in considering dynamics in Banaach spaces, we give up unique notions of orthogonality and Riesz dual pairings (we discuss contraction rates in uniformly convex spaces admitting generalized dual pairings in \ref{sec:regression}).

\begin{definition}[Cauchy problem on a Banach space]
Let $V$ be a Banach space and $f(t, u) : \R^+ \times V \to V$ a nonlinear time-varying operator on $V$. The \textit{Cauchy problem} for $f$ is to identify, for some $u_0 \in V$, some $u(t) : \R^+ \to V$ such that:
\begin{equation}
\begin{split}
    \label{eq:cauchy}
    \dot{u} &= f(t, u)\ \forall t \ge 0\\
    u(0) &= u_0
\end{split}
\end{equation}
\end{definition}

We refer the reader to \cite{ichikawa1979existence} for a reference on the well-posedness of the Cauchy problem for nonlinear evolution equations in Banach spaces.
We assume the existence of particular classical ($t$-differentiable) solutions $u_1(t), u_2(t)$ to \eqref{eq:cauchy} as needed, and that the set $E \subset V$ on which \eqref{eq:cauchy} is well-posed is a closed subspace of $V$, subsequently using $V$ synonymously with $E$.

For a Cauchy problem in a Banach space $V$, the contraction condition \eqref{eq:contracting_sys} corresponds to satisfy an exponential norm bound on particular solutions:
\begin{align}
    \label{eq:contracting_family}
    \norm{u_1(t) - u_2(t)} \le e^{\lambda(t-s)}\norm{u_1(s) - u_2(s)}
\end{align}
for all $t \ge s \ge 0$ and some $\lambda \in \R$. We introduce some definitions which enable the verification of \eqref{eq:contracting_family} in Banach spaces.

\begin{definition}[Gateaux semi-inner product on a Banach space]
Let $V$ be a real Banach space; the right and left \textit{semi-inner products} (in the sense of \cite{lumer1961semi}) are defined as
\begin{align}
\label{eq:sip}
(u, v)_{\pm} = \norm{u} \lim_{h\to0^{\pm}}\frac{\norm{u + hv} - \norm{u}}{h}
\end{align}
\end{definition}
Definition \eqref{eq:sip} is sometimes referred to as the \textit{Gateaux formula} (see \cite{davydov2021non}), due to the inner limit which is the left- or right-Gateaux derivative of the norm $\norm{\cdot}_V$, and exists in any Banach space (but may not be equal). Thus $(u, v)_+ = (u, v)_-$ when $\norm{\cdot}_V$ is Gateaux-differentiable, e.g. $V = L^p$ for $1 < p < \infty$, and $(u, v)_+ = (u, v)_- = \langle u, v\rangle$ in $L^2$. We refer the reader to \cite{lumer1961semi} for properties, with the notable ones being positive-definiteness, norm-compatibility, and Cauchy-Schwarz while losing symmetry and bilinearity in general. 

We employ the following key result (see variations in \cite{soderlind2006logarithmic, aminzare2013logarithmic, aminzare2014contraction}) which is equivalent to the contraction condition \eqref{eq:contracting_family}:
\begin{theorem}[Strongly negative-definite semi-inner product implies contraction]
Let $V$ be a real Banach space and $f \in C^1(V)$; define the functional
\begin{align}
    \label{eq:llc}
    M(f) &:= \sup_{u \ne v \in V} \frac{(u-v, f(u) - f(v))_+}{\norm{u-v}^2}
\end{align}
If an evolution family $\phi(s, t, u)$ is generated by a differential equation $\dot{u} = f(t, u)$, then for any $u, v \in V$ and $s \le t \in \R^+$, 
\begin{align*}
    \norm{\phi(s, t, u) - \phi(s, t, v)} \le e^{\int_s^t M(f_{t'})dt'}\norm{u - v}
\end{align*}
\begin{proof}
First, we apply the Dini identity (see 5.3, \cite{soderlind2006logarithmic})
\begin{align}
\label{eq:dini_bound}
D^+_t\norm{u} = \frac{(u, \dot{u})_+}{\norm{u}^2}\norm{u}
\end{align}
where $D^+$ denotes the upper Dini derivative (which exists for any norm, but may be infinite)
\begin{align}
\label{eq:dini}
D^+_t \varphi(t) := \limsup_{h\to 0^+} \frac{\varphi(t + h) - \varphi(t)}{h}
\end{align}
along with definition \eqref{eq:llc} to obtain, for any two solutions $u(t), v(t)$,
$$
D^+_t \norm{u(t)-v(t)} \le M(f_t)\norm{u(t)-v(t)}
$$
where $f_t(u) = f(t, u)$. By applying a Gr\"onwall inequality for Dini derivatives (see Lemma 11, \cite{davydov2021non}), we have:
\begin{align}
    \label{eq:llc_contraction}
    \norm{u(t) - v(t)} \le e^{\int_s^tM(f_{t'})dt'}\norm{u(s) - v(s)}
\end{align}
and the result follows.
\end{proof}
\end{theorem}
In the special case where $\phi(s, t, u) := \Phi(t)u$ where $\Phi(t)$ is a $C^0$-semigroup, we have that $M(f) < 0$ implies the existence of $t_0$ such that $\norm{\Phi(t_0)} < 1$, implying that $\Phi(t)$ is a contraction semigroup.  

Classical results in contraction theory by \cite{lohmiller1998contraction} as well as numerical and functional analysis by \cite{dahlquist1958stability, ladas1972differential} characterize \eqref{eq:llc} instead in a \textit{differential} sense, by estimating contraction rates using the Fr\'echet linearization of the vector field at each point.
\begin{definition}[Fr\'echet derivative]
A continuous map $f: V \to W$ between Banach space $V, W$ is \textit{Fr\'echet-differentiable} at a point $u \in V$ if there exists a bounded linear operator $Df(u) \in \B(V, W)$ satisfying:
\begin{align*}
    \lim_{\norm{h}\to 0}\frac{\norm{f(u+h) - f(u) - Df(u)h}}{\norm{h}} = 0
\end{align*}
Regularity classes $f \in C^k(V)$ are defined accordingly and partial derivatives $\frac{\partial f}{\partial u}(u, v, ...)$ are defined as Fr\'echet derivatives with other arguments supplied.
\end{definition}
By an application of the fundamental theorem of calculus for Fr\'echet derivatives, it is sufficient (see \cite{ladas1972differential}, Chapter 5.4  or Appendix \ref{thm:orig_contraction} for proof) to show that the linearization is contracting at all points, i.e.
\begin{align}
    \label{eq:linearized_contraction}
    \sup_{u \in V}M(Df(u)) = \sup_{u, v \in V}\frac{(v, Df(u)v)_+}{\norm{v}^2} = \lambda < 0
\end{align}
in order to obtain the same growth bound between solutions as \eqref{eq:llc_contraction}. (We refer to \eqref{eq:linearized_contraction} as the \textit{differential} contraction rate.) This bound is looser, but arises from the fact that the operator norm of $Df(u)$ is upper-bounded by the Lipschitz constant of $f$, which in turn can be used (via the \textit{logarithmic Lipschitz constant} as in \cite{aminzare2013logarithmic}) to obtain an upper bound for \eqref{eq:llc}.

This is the foundation of contraction analysis as often applied in control \& stability theory in $\R^n$ and $\C^n$, since there the condition \eqref{eq:linearized_contraction} reduces to the linear matrix inequality:
\begin{align}
    \label{eq:contraction_rn}
    \sup_{u, v \in V} \frac{\Re\langle v, J_f(u)v\rangle}{\langle v, v\rangle} &= \sup_{u \in V}\max\left(\Spec\left(\frac{J_f(u) + J_f(u)^*}{2}\right)\right) < 0
\end{align}
where $J_f(u)$ is the Jacobian of $f$ at $u$, $\Spec(A)$ is the set of eigenvalues of a matrix $A$, and $A^*$ is the adjoint or conjugate-transpose. The quantity $\langle v, J_f(u)v\rangle$ in \eqref{eq:contraction_rn} is better known as the \textit{numerical range} of $J_f(u)$. In Hilbert spaces, this inequality is related to the notion of a \textit{dissipative} operator.

\begin{definition}[Dissipative operator in a Hilbert space]
An operator $A : H \to H$ with domain $D(A) \subset H$ dense in a Hilbert space $H$ is \textit{dissipative} if 
\begin{align}
    \label{eq:dissipative}
    \Re\langle v, Av\rangle \le 0\quad \forall x \in H
\end{align}
\end{definition}
By this definition, $J_f(u)$ in \eqref{eq:contraction_rn} is semi-contractive ($M(J_f(u)) \le 0$) if and only if $J_f(u)$ is dissipative. 
A well-known example in the continuum setting (see \cite{soderlind2006logarithmic}) is the Laplacian $\Delta$, which is dissipative on the set of compactly supported smooth functions $C_c^\infty(\Omega) \subset H_0^2(\Omega)$, a fact which we apply in the analysis of diffusive nonlinear PDEs and feedback connections. Definition \eqref{eq:dissipative} and the semi-contraction equivalence generalizes straightforwardly to Banach spaces via substitution of the semi-inner product $(v, Av)_+$ and application of the Dini identity \eqref{eq:dini_bound}. Thus, the weighted semi-inner products we introduce in the Section \ref{sec:weighted_contraction} generalize this notion of a dissipative operator (and of dissipative dynamics).

In an arbitrary Banach space $V$, the contraction rate of a bounded linear operator \eqref{eq:linearized_contraction} is equal to the following functional (see Theorem 1, \cite{aminzare2013logarithmic}):
\begin{align}
    \label{eq:lognorm}
    \mu(A) &:= \lim_{h \to 0^+}\frac1h\left(\norm{I + hA})_{\B(V)} - 1\right)
\end{align}
which is well-known as the \textit{logarithmic norm} or \textit{matrix measure} (see \cite{dahlquist1958stability, soderlind2006logarithmic}). This equality is non-trivial to show; we mention this definition for historical reasons, as it arises naturally when taking the upper Dini derivative of the norm of a perturbation (see Appendix \ref{thm:orig_contraction}). However, the approach we present in this paper is developed from the semi-inner product \eqref{eq:llc} (henceforth referred to as SIP-contracting or simply contracting) rather than the logarithmic norm \eqref{eq:lognorm} for the following reasons. 

First, estimating the differential contraction rate \eqref{eq:linearized_contraction} becomes problematic in the infinite-dimensional case for unbounded linear operators (or nonlinear operators whose linearization is unbounded) such as differential operators, since we assume that $\norm{I + hA} < \infty$. The Hilbert space formula \eqref{eq:contraction_rn} allows, for instance, evaluating the contractivity of the Laplacian using integration by parts (see \cite{soderlind2006logarithmic, cisneros2020contraction}), and the semi-inner product \eqref{eq:sip} definition of the numerical range extends these quadratic forms to unbounded operators on Banach spaces.

Secondly, the differential contraction rate \eqref{eq:linearized_contraction} is looser than the integral contraction rate \eqref{eq:llc_contraction}, since $M(f) \le \sup_u M(Df(u))$ (see \cite{aminzare2014contraction}). In particular, converse contraction results hold for the latter, i.e. \eqref{eq:contracting_family} implies that $M(f_t) \le \lambda$ (see Proposition 3, \cite{aminzare2014contraction}) but no such result holds for $\sup_u M(Df_t(u))$. We still use differentiation when convenient (e.g. when analyzing Jacobians of feedback systems), but otherwise use $M(f_t)$. The integral contraction rate \eqref{eq:llc_contraction} is also tighter than the \textit{logarithmic Lipschitz constant}
\begin{align*}
    LLC(f) &:= \lim_{h \to 0^+} \frac{L(I + hf) - 1}{h} \ge M(f)
\end{align*}
which may be seen as a natural extension of the logarithmic norm to nonlinear arguments (see \cite{aminzare2014contraction}). However, \eqref{eq:llc_contraction} does not require Lipschitz continuity; only a \textit{one-sided Lipschitz}-like condition (see \cite{davydov2021non} for further discussion of this requirement and its relationship to others such as the Demidovich condition and operator measures).

Lastly, the SIP constant \eqref{eq:llc} carries over to normed spaces several properties of the numerical range of operators in Hilbert spaces, such as subadditivity and homogeneity for positive scalar multipliers in the second argument (see \cite{lumer1961semi}), which we rely on for the analysis of coupled systems (Section \ref{sec:combinations}).

\subsection{PDEs as ODEs on Banach spaces}
\label{sec:ode_pde}
While our exposition is from the perspective of abstract differential equations in Banach spaces, we briefly discuss its application to partial differential equations. In general, converting a PDE of the form
\begin{align}
    \label{eq:pde0}
    \partial_t u &= f(t, \partial^{\alpha^1} u, ..., \partial^{\alpha^n} u)
\end{align}
where $\alpha^k$ are multi-indices to a well-posed Cauchy problem \eqref{eq:cauchy} on a Banach space $V$, where $F$ may be a nonlinear (differential) operator, is not possible when $D(F) \ne V$. Moreover, spaces of (classically) differentiable functions are typically not complete. General procedures for obtaining suitable ``extensions'' of \eqref{eq:pde0} to an appropriate Banach space $V$ is outside the scope of our work; however, we mention some common approaches which we apply in the analysis of PDEs (Section \ref{sec:pde}). For bounded linear $A$ on some domain $D(A) \subset V$ which is dense, we can study \eqref{eq:cauchy} with its unique continuous linear extension (for example, the extension of the Laplacian $\Delta$ from  $C_c^\infty$ to the set of trace-zero functions $H_0^1$). Contraction results of the ``extended'' equation then apply to solutions of the original equation, assuming existence. More generally, the Hille-Yosida theorem (see \cite{pazy2012semigroups}) gives necessary and sufficient conditions for linear $A$ with $D(A)$ a dense subspace of $V$ to generate a $C^0$-semigroup of operators on $V$; extensions to nonlinear semigroups are developed in \cite{ladas1972differential}. A second, related approach is to relax the notion of a solution $u$ from a function possessing sufficiently many (classical) space derivatives to one possessing sufficiently many \textit{weak} derivatives in the sense of test functions; this weak formulation of \eqref{eq:pde0} may admit a representation as an ODE in a Sobolev space $W^{k, p}$ subject to boundary constraints. We note that an advantage of the semi-inner product formulation of the contraction rate \eqref{eq:llc} is that we do not necessarily require bounded or Lipschitz extensions $\overbar{F}$ of $F$ to the domain $V$, only that $\overbar{F}$ is one-sided Lipschitz or that $D\overbar{F}(u)$ has uniformly bounded real numerical range.

\subsection{Notation}
We use $V,W$ to denote Banach spaces unless otherwise mentioned. The set $\B(V)$ denotes the Banach space of bounded linear operators $A: V \to V$ with the standard operator norm, and $\B(V, W)$ if the domain and range are $V \ne W$.
We use the following shorthands for partial derivatives of time-varying functions on $V$:
\begin{align}
    \label{eq:nonautonomous_shorthand}
    f_t(u) &:= f(t, u),\quad \dot{f}_t(u) := \frac{\partial f}{\partial t}(t, u),\quad Df_t(u) := \frac{\partial f}{\partial u}(t, u)
\end{align}
where $D$ is the Fr\'echet derivative. Finally, we refer to $M(f)$ in \eqref{eq:llc} simply as the \textit{contraction rate} of $f$ throughout (where \textit{contraction} by convention corresponds to $M(f) < 0$).

\section{Contraction in weighted spaces}
\label{sec:weighted_contraction}

One of the important extensions to classical contraction theory is the relaxation of contraction in the $\ell^2/L^2$ norm (referred to in \cite{lohmiller1998contraction} as the \textit{identity metric}) to contraction in geodesic distances induced by time-varying Riemannian metrics. The main result is that differential (in the sense of \eqref{eq:linearized_contraction}) contraction in a Riemannian metric (given by a smoothly-weighted inner-product) implies contraction in the geodesic distance:
\begin{equation}
\begin{split}
    \label{eq:contraction_metrics}
    \langle v, w\rangle_{T_u\R^n} &:= \langle \Theta(t, u)v, \Theta(t, u)w\rangle_{\R^n}\\
    \frac{d}{dt}\langle \delta u, \delta u\rangle_{T_u\R^n} &\le \lambda \langle \delta u, \delta u\rangle_{T_u\R^n} \implies d_g(u(t), v(t)) \le Ce^{\lambda t}d_g(u(0), v(0))
\end{split}
\end{equation}
(See \cite{lohmiller1998contraction} and Theorem 2.3,  \cite{simpson2014contraction}).

While this construction necessarily assumes the model space is Hilbertian, we propose several generalizations to normed spaces using \textit{weighted} semi-inner products applied to integral \eqref{eq:llc} as well as differential \eqref{eq:linearized_contraction}, \eqref{eq:contraction_metrics} contraction rates.

\subsection{Contraction under linearly weighted semi-inner products}
\begin{theorem}[Contraction in a constant weight]
Let $\Theta \in GL(V, W)$ be a bijective bounded linear operator between Banach spaces $V, W$ and define the $\Theta$-weighted semi-inner product:
\begin{align}
    \label{eq:integral_weight}
    (u, v)_+^{\Theta} := (\Theta u, \Theta v)_+
\end{align}
Define $M^{\Theta}$ to be the contraction constant \eqref{eq:llc} with respect to this weighted SIP. Then:
\begin{enumerate}
    \item For any $f : V \to W$, we have $M^\Theta(f) = M(\Theta\circ f \circ \Theta^{-1})$
    \item For solutions $u_1(t), u_2(t)$ to $\dot{u} = f(t, u)$, we have:
    $$
    \norm{u_1(t) - u_2(t)} \le \kappa(\Theta)e^{M^{\Theta}(f_t)}\norm{u_1(0)-u_2(0)}
    $$
    where $\kappa$ is the condition number of $\Theta$.
\end{enumerate}
\end{theorem}
\begin{proof}
By direct substitution:
\begin{align}
    \label{eq:conj_llc}
    M^{\Theta}(f) &= \sup_{u\ne v \in W}\frac{(u - v, \Theta f(\Theta ^{-1}u) - \Theta f(\Theta^{-1}v))_+}{\norm{u-v}^2} = M(\Theta \circ f \circ \Theta^{-1})
\end{align} 
Next, suppose the weighted contraction rate of $f_t$ satisfies $M^\Theta(f_t) \le \lambda$. By linearity and \eqref{eq:conj_llc}, this is equivalent to showing that the conjugate dynamics $v = \Theta u$ is contracting in $W$. Hence for any two solutions $u_1(t), u_2(t)$:
\begin{align*}
    \norm{u_1(t) - u_2(t)} &\le \norm{\Theta^{-1}}\norm{\Theta u_1(t) - \Theta u_2(t)} \le \norm{\Theta^{-1}}\norm{\Theta} e^{\lambda t}\norm{u_1(0) - u_2(0)}
\end{align*}
\end{proof}
Indeed, the same argument applies in the reverse direction; suppose that $M(f_t) \le \lambda$, then:
\begin{align*}
    \norm{\Theta u_1(t) - \Theta u_2(t)} &= \kappa(\Theta) e^{\lambda t} \norm{\Theta u_1(0) - \Theta u_2(0)}
\end{align*}

This can be used to establish the following coordinate independence of asymptotic contraction rates. For any weights $\Theta_1, \Theta_2 \in GL(V, W)$, we have:
\begin{equation}
\begin{split}
    \label{eq:so_equivalence}
    &\norm{u_1(t) - u_2(t)}_{\Theta_1} \le \kappa(\Theta_1)\kappa(\Theta_2)e^{\lambda_{\Theta_2} t}\norm{u_1(0) - u_2(0)}_{\Theta_1}
\end{split}
\end{equation}
where $\lambda_{\Theta_2}$ is the contraction rate with respect to weight $\Theta_2$. Thus, a flow is contracting after a \textit{finite overshoot} (adapting terminology of \cite{sontag2014three}) in a weighted SIP if and only if it is contracting after a finite overshoot in \textit{any} weighted SIP. 

\subsubsection{Metric dependence}
\label{sec:metric_dependence}
Statement \eqref{eq:so_equivalence} establishes the equivalence of asymptotic, rather than time-uniform, exponential convergence of any two solutions under any coordinate change $v = \Theta u$. This result is analogous to the invariance of Lyapunov exponents under measure-preserving invertible transformations as in Oseledets' Theorem (see \cite{ruelle1979ergodic}); both the contraction rate after a finite overshoot and Lyapunov exponents capture asymptotic properties of the dynamics.

On the other hand, it is known that the time-dependent contraction rate $M(f_t)$ in general depends on the choice of norm (see 3.5-3.6, \cite{banasiak2020logarithmic}). For completeness, we discuss how a naive comparison argument between $M^\Theta(f_t)$ and $M(f_t)$ fails.

Since $\Theta$ is a bijection $V \to W$, we have the two facts
\begin{align}
\label{eq:inverse_bounds}
\inf_{v\ne 0 \in V} \frac{\norm{\Theta v}}{\norm{v}} = \frac{1}{\norm{\Theta^{-1}}} \text{ and } \norm{\Theta}\norm{\Theta^{-1}} = \kappa(\Theta) \ge 1
\end{align}
We can then deduce that
\begin{align*}
    M^\Theta(f) &= \sup_{u\ne v\in V}\lim_{h\to 0^+}\frac1h\left(\frac{\norm{\Theta(u-v+hf(u)-hf(v))}}{\norm{\Theta(u-v)}}-1\right) \\
    &\le  \kappa(\Theta)M(f) + \lim_{h \to 0}\frac{\kappa(\Theta) - 1}{h} \\
    M^\Theta(f) &\ge \kappa(\Theta)^{-1}M(f) + \lim_{h \to 0}\frac{\kappa(\Theta)^{-1} - 1}{h}
\end{align*}
However, from \eqref{eq:inverse_bounds} we have that $\kappa(\Theta) \ge 1$ with equality iff $\Theta \propto I_{VW}$ where $I_{VW}$ is an isometric isomorphism. We can deduce if so that $M^\Theta(f) = M(f)$, and otherwise only that $-\infty \le M^\Theta(f) \le \infty$.

\subsection{Contraction under Lipschitz coordinate transformations}
We now generalize \eqref{eq:conj_llc} to the nonlinear time-varying setting.

\begin{theorem}[Nonlinear coordinate change]
Let $\theta(t, u) \in C^1(\R^+ \times V,  W)$ be a smoothly time-varying family of maps between Banach spaces $V, W$ such that:
\begin{enumerate}
    \item $\theta_t$ is bijective for all $t$
    \item $\theta_t, \theta_t^{-1}$ have $t$-uniformly bounded Lipschitz constants $L(\theta_t), L(\theta_t^{-1})$
\end{enumerate}
Let $\dot{u} = f(t, u)$ and consider the nonlinear coordinate transformation $v = \theta(t, u)$. If $v(t)$ is contracting with rate $\lambda$, then $u(t)$ is also contracting with rate $\lambda$ after a finite overshoot.
\end{theorem}
\begin{proof}
The conjugate dynamics of $v$ are:
\begin{align}
\label{eq:conj_nonlinear}
\dot{v} &= g(t, v) = \dot{\theta}_t(\theta_t^{-1}(v)) + D\theta_t(\theta_t^{-1}(v))f_t(\theta_t^{-1}(v))
\end{align}
By hypothesis, we have $\sup_t M(g_t) = \lambda$; then for any two solutions $u_1(t), u_2(t)$ we have:
\begin{align*}
    \norm{u_1(t) - u_2(t)}_V &\le L(\theta_t^{-1})\norm{v_1(t) - v_2(t)}_W \le L(\theta_t^{-1})L(\theta_0) e^{\lambda t} \norm{u_1(0) - u_2(0)}_V 
\end{align*}
where $L(\Theta_t^{-1})$ is uniformly bounded by assumption.
for some $C > 0$. 
\end{proof}
Note that in general, \eqref{eq:conj_nonlinear} does not arise from a weighted SIP. Furthermore, the converse result follows from $t$-uniform boundedness of $L(\Theta_t)$, resulting in an identical finite-overshoot contraction equivalence as in the case of a linear autonomous coordinate change \eqref{eq:so_equivalence}.

\subsection{Contraction under locally linearly weighted SIPs}
\label{sec:local_linear_latent}
We now give a generalization of geodesic contraction in Riemannian metrics to Banach spaces by emulating duality on the tangent space with semi-inner products. Formally, we consider $V$ to be a Banach manifold with $T_pV \isoto V$ via the trivial isomorphism. We imbue $T_pV$ with a smoothly varying semi-inner product space structure via the weighted SIP:
\begin{align}
    \label{eq:differential_weight}
    (u, v)_+^{\Theta} &:= (\Theta(t, p)u, \Theta(t, p)v)_+
\end{align}
which generalizes the prior weighted SIP \eqref{eq:integral_weight} in that ``curvature'' is no longer zero. 

(We use terms such as ``metric'' and ``curvature'' informally; there are several possible notions of angle and geodesic distance -- for instance, via symmetrized Gateaux derivatives (Section 3.6, \cite{balestro2017angles}) -- one might define in a semi-inner product space, and no presently known canonical choice, unlike in Euclidean-like spaces. We do not make use of these distances explicitly, opting instead to measure the contraction rate \eqref{eq:differential_weight} via the differential \eqref{eq:linearized_contraction} rather than integral \eqref{eq:llc} condition in the following theorem. As we show below, differential contraction implies the asymptotic convergence of solutions in any geodesic-like distance which preserves the topology of $V$.) 

\begin{theorem}[Differential contraction in a smoothly-varying weight]
As above, let $(\cdot, \cdot)_+^\Theta$ be a weighted semi-inner product with $\Theta(t, u)$ such that:
\begin{enumerate}
    \item $\Theta(t, u) \in C^1(\R^+\times V, GL(V, W))$ for some Banach space $W$
    \item $\sup_{t, u}\max(\norm{\Theta(t, u)}, \norm{\Theta(t, u)^{-1}}) = C < \infty$
\end{enumerate}
Then if $\dot{u} = f(t, u)$ and $M^\Theta(f_t(u)) \le \lambda$ then for any two solutions $u_1(t), u_2(t)$, we have:
\begin{align*}
    \norm{u_1(t) - u_2(t)} \le C^2e^{\lambda t}\norm{u_1(0) - u_2(0)}
\end{align*}
\end{theorem}
\begin{proof}
The weighted SIP \eqref{eq:differential_weight} defines a parametric norm $\norm{\cdot}_{t, u}$ on $W$. We start by applying the Dini derivative identity \eqref{eq:dini_bound} in this norm:
\begin{align}
    \label{eq:weighted_dini_bound}
    D_t^+\norm{\delta u}_{t,u} &= \frac{\left(\Theta(t, u)\delta u, \frac{d}{dt}[\Theta(t, u)\delta u]\right)_+}{\norm{\delta u}_{t,u}^2} \norm{\delta u}_{t,u}
\end{align}
where $\delta u \in T_uV$ is a perturbation. Let $\frac{d}{dt}[\Theta(t, u)\delta u = G(t,u)\delta u$, where
\begin{align*}
    G_t(u) &= \dot{\Theta}_t(u) + \Theta_t(u) Df_t(u) \in \B(V, W)
\end{align*}
Letting $M^{\Theta}$ be the weighted contraction rate induced by \eqref{eq:differential_weight}, we have the equality 
\begin{align*}
    M^{\Theta}(\Theta(t,u)^{-1}G(t, u)) &= M(G(t, u)\Theta(t,u)^{-1})
\end{align*}
By hypothesis, this contraction rate is upper-bounded by $\lambda$. Applying definition \eqref{eq:llc}, we have that \eqref{eq:weighted_dini_bound} reduces to:
\begin{align*}
    D_t^+\norm{\delta u}_{t,u} &\le \lambda \norm{\delta u}_{t,u} 
\end{align*}
Applying Gr\"onwall's inequality, the perturbation in the identity metric then satisfies the growth bound:
\begin{align*}
    \norm{\delta u(t)}_{t, u(t)} &\le e^{\lambda t}\norm{\delta u(0)}_{0, u(0)} \implies \norm{\delta u(t)} \le C^2e^{\lambda t}\norm{\delta u(0)}
\end{align*}
Applying the fundamental theorem of calculus for Fr\'echet spaces as in  \cite{ladas1972differential} (5.4), we can conclude that any two solutions $u_1(t), u_2(t)$ are contracting with rate $\lambda$ after a finite overshoot.
\end{proof}

\subsection{Contraction in latent complex-valued spaces}
\label{sec:complex_contraction}
Many systems are best described in complex vector spaces, either intrinsically or by complexification. Prior work on logarithmic norms and semi-inner products for incremental stability analysis (e.g. \cite{aminzare2013logarithmic, aminzare2014contraction, soderlind2006logarithmic, ladas1972differential}) generally assumes the underlying field is real, which makes the right- and left-Gateaux differentials in the expression: 
\begin{align}
    \label{eq:right_sip}
    (u, v)_+ = \norm{u}\lim_{h\to 0^+}\frac{\norm{u + hv} - \norm{u}}{h}
\end{align}
well-defined.
However, the extension of this Gateaux formula SIP to complex Banach spaces, in the sense of Gateaux-holomorphic functionals, is problematic since the norm is not complex-(Gateaux-) differentiable (e.g. $|z|^2, z \in \C$). 

We list some solutions to this issue. First, we can simply choose a Banach space $V$ whose norm $\norm{\cdot}_V$ is Gateaux-differentiable (the left- and right-limits in \eqref{eq:right_sip} are equal), e.g. $L^p$ for $1 < p < \infty$. Then,
\begin{align}
    (u, v)_{V(\C)} &:= (u, v)_{V(\R)} + i(u, iv)_{V(\R)}
\end{align}
induces a true semi-inner product on the complex space $V(\C)$ in the sense of complex-homogeneity (see Proposition 28, \cite{zhang2009reproducing} for proof). In this case, we have that 
\begin{align*}
    M_{V(\C)}(f) &= \sup_{u\ne v}\frac{\Re (u-v, f(u)-f(v))_{V(\C)}}{\norm{u-v}^2} = \sup_{u\ne v}\frac{(u-v, f(u)-f(v))_{V(\R)}}{\norm{u-v}^2}\\
    &= M_{V(\R)}(f)
\end{align*}
thus the contraction rate is unchanged.

If the norm is not $G$-differentiable (e.g. $L^1$), we may invoke a comparison argument with contraction in a complex Hilbert space (e.g. $L^2(\C^n)$) as we describe in Section \ref{sec:lp_estimates}, and measure the contraction rate of the linearization, which reduces to the $L^2$-logarithmic norm:
\begin{align}
    \label{eq:l2_lognorm}
    M_2(Df(u)) &= \sup_{v \in V}\frac{\Re \langle v, Df(u)v\rangle}{\langle v, v\rangle}
\end{align}
(If additionally \eqref{eq:l2_lognorm} is weighted by $\Theta(t, u)$ as in \eqref{eq:differential_weight}, then 
$$
M(t, u) = \Theta(t, u)^*\Theta(t, u)
$$ is the Riemannian metric on $L^2$ with respect to which geodesic contraction is guaranteed.)

Finally, we can measure the contraction rate of the linearization at all points, but split into two operators representing the real and complex parts:
\begin{align*}
    \frac{d}{dt}\begin{bmatrix}\delta u_r \\ \delta u_i\end{bmatrix} &= \begin{bmatrix}F_r & -F_i \\ F_i & F_r\end{bmatrix}\begin{bmatrix}\delta u_r \\ \delta u_i\end{bmatrix}
\end{align*}
where $F_r$ is a shorthand for $\Re[Df_t(u)]$. We then apply a suitable norm on the product space $\Re[V]^2$ and measure the contraction rate in this norm, either directly or by using the feedback combination properties (section \ref{sec:feedback_combinations}). Contraction in this product norm implies contraction in the complex space $V$, but the converse may not be true for general norms on $V$ and $\Re[V]^2$.

\subsection{$L^p$ estimates from $L^2$ contraction}
\label{sec:lp_estimates}
We now show some cases in which the $L^p$ contraction rate estimated in $L^2$, which is often easier to calculate by the simple formula \eqref{eq:l2_lognorm} on the linearization.

Let $E \subset \R^n$ be a finite-measure set (we use $m$ to denote the Lebesgue outer measure) and $\dot{u} = f(t, u)$ an ODE on $V = L^p(E)$. Let $\delta u(t) = T_{u(t)}V$ be any perturbation of a solution. Throughout, we assume the following $L^2$ contraction rate to be finite:
\begin{align}
    \label{eq:l2_bound}
    \lambda_2 &:= \sup_t M_2(f_t)
\end{align}

\textbf{Case $1 \le p < 2$.} By H\"older's inequality, we have
\begin{align}
    \norm{\delta u(t)}_p &\le m(E)^{\frac1p - \frac12}\norm{\delta u(t)}_2 \le m(E)^{\frac1p - \frac12}e^{\lambda_2 t}\norm{\delta u(0)}_2
\end{align}

Thus, for initial conditions $u_1, u_2 \in L^p \cap L^2$, we can bound their $p$-distance at $t$ as a constant multiple of their initial $2$-distance. This does not give us global contraction in $L^p$ since in general, $\norm{\delta u(0)}_2$ may be infinite despite finite $L^p$ initial conditions (e.g. $p = 1$, $E=[0,1]$, $\delta u(0, x) = x^{-1/2}$). 

However, we can make this result more precise in the Banach space $L^p(E) \cap L^2(E)$ with the norm $\norm{\cdot}_{p, 2} = \max(\norm{\cdot}_p, \norm{\cdot}_2)$:
\begin{align*}
    \norm{\delta u(t)}_{p,2} \le \max(1, m(E)^{\frac1p - \frac12})e^{\lambda_2 t}\norm{\delta u(0)}_2 \le \max(1, m(E)^{\frac1p - \frac12})e^{\lambda_2 t}\norm{\delta u(0)}_{p, 2}
\end{align*}
We can also extend this to contraction in a weighted inner-product (see section \ref{sec:local_linear_latent}) as follows. Let $\Theta(t, u) \in C^1(\R^+\times L^2(E), GL(L^2(E)))$ and suppose the contraction rate $\lambda_2^\Theta$ is with respect to this weight. Define $C_\Theta = \norm{\Theta(0, u(0))}\sup_t\norm{\Theta(t, u(t))}^{-1}$. Then by the same argument as in Section \ref{sec:local_linear_latent},
\begin{align*}
    \norm{\delta u(t)}_{p,2} &\le C_\Theta\max(1, m(E)^{\frac1p - \frac12})e^{\lambda_2^\Theta t}\norm{\delta u(0)}_{p,2} 
\end{align*}

\textbf{Case $2 < p \le \infty$}. Here we make the assumption that there exists an $L^\infty$ bounding box for $u(t)$, i.e. $\sup_t \norm{\delta u(t)}_\infty = B < \infty$. By H\"older again, 
\begin{align*}
    \norm{\delta u(t)}_p &\le \norm{\delta u(t)}_\infty^{1-\frac{2}{p}}\norm{\delta u(t)}_2^{\frac{2}{p}}  \le B^{1-\frac{2}{p}} \left(e^{\lambda_2 t}\norm{\delta u(0)}_2\right)^{\frac{2}{p}} \\
    &\le B^{1-\frac{2}{p}} \left(e^{\lambda_2 t}m(E)^{\frac12 - \frac1p}\norm{\delta u(0)}_p\right)^{\frac{2}{p}}
\end{align*}
We can also extend this to a weighted space, with $\Theta, \lambda_2^\Theta, C_\Theta$ as before:
\begin{align*}
    \norm{\delta u(t)}_p &\le B^{1-\frac{2}{p}} \left(C_\Theta e^{\lambda_2 t}m(E)^{\frac12 - \frac1p}\norm{\delta u(0)}_p\right)^{\frac{2}{p}}
\end{align*}

\section{Contraction to invariant subsets}
We now apply contraction in weighted semi-inner products to derive several notions of contraction to subsets (vector subspaces as in \cite{pham2007stable} and submanifolds as in \cite{manchester2017control}) for dynamics in Banach spaces. 

\subsection{Contraction to a subspace}
\label{sec:subspace_contraction}
\begin{theorem}[Contraction to subspaces via weighted semi-inner products]
Let $V$ be a Banach space, and $P$ a bounded linear projection operator on $V$ (i.e. $P^2 = P$). Let $Q = I - P$, and $M^Q(A) = \sup_{v \in V}\frac{(Qv, QAv)_+}{\norm{Qv}^2}$ be the $Q$-weighted contraction rate. Let $\dot{u} = f(t, u)$ and suppose that:
\begin{enumerate}
    \item $\im(P)$ is a flow-invariant subspace of $u(t)$
    \item $\sup_{t, u \in V}M^Q(Df_t(u)) = \lambda 0$
\end{enumerate}
Then, $u(t)$ is contracting to the subspace $\im(P)$ in the set distance $d(u, \im(P)) = \inf_{v\in \im(P)}\norm{u - v}$.
\end{theorem}
\begin{proof}
$P$ gives a decomposition of $V$ into a direct sum of closed subspaces as $V = \im(P) \oplus \ker(P)$; moreover, $Q$ is also a projection, and satisfies $\im(Q) = \ker(P)$. Consider the projected system $v = Qu$ with dynamics 
\begin{align}
    \label{eq:subspace_dynamics}
    \dot{v} &= Qf(t, u) = Qf(t, v + Pu)
\end{align}
By hypothesis:
\begin{enumerate}
    \item Since $\im(P)$ is a flow-invariant subspace of $u(t)$, we have
    \begin{align}
        \label{eq:invariant_subspace}
        Qf(t, Pv) = 0\ \forall t, v \in V
    \end{align}
    thus $v(t)$ has a fixed point at $v = 0$. 
    \item Since $u(t)$ is contracting in the $Q$-weighted SIP, we have
    \begin{equation}
    \begin{split}
        \label{eq:subspace_contraction}
        \sup_{t, u \in V}M^Q(Df_t(u))  = \sup_{t, u, w \in V} \frac{\left(Qw, QDf_t(u)w\right)_+}{\norm{Qw}^2} < 0
    \end{split}
    \end{equation}
    Then by linearity, $\delta v = Q\delta u$, and $v(t)$ is contracting; that is, a perturbation $\delta v$ with dynamics
    $$
    \frac{d}{dt}\delta v = Q Df_t(u)\delta v
    $$
    is exponentially stable.
\end{enumerate}
Thus, all solutions $v(t)$ of \eqref{eq:subspace_dynamics} approach the fixed point $0$ exponentially. Then, $u(t)$ is contracting to the subspace $\im(P)$.
\end{proof} 
Note that in the $L^2$ case, if $Q = A^*A$ such that $AA^* = I$, then \eqref{eq:subspace_contraction} reduces to: 
\begin{align*}
    M^Q(Df_t(u)) &= \sup_{w \in V} \frac{\Re \langle Qw, QDf_t(u)w\rangle}{\langle Qw, Qw\rangle} = \sup \Re[\Spec(ADf_t(u)A^*)]
\end{align*}
which reduces to the prior known condition in $\R^n$ given by \cite{pham2007stable} $ADf_t(u)A^\intercal < 0$. On the other hand, \eqref{eq:subspace_contraction} does not require a notion of orthogonality or self-duality.

We can further relax condition \eqref{eq:subspace_contraction} to contraction in a latent space. Let $\Theta(t, v) : C^1(\R^+\times V, GL(V))$ be a uniformly bounded family of operators as in \ref{sec:local_linear_latent}; if the $\Theta$-weighted linearized $v$ dynamics are contracting in the $Q$-weighted SIP,
\begin{align}
    \label{eq:weighted_subspace_contraction}
    \sup_{t, u\in V} M^{Q}((\dot{\Theta}_t(u) + D\Theta_t(u)Df_t(u))\Theta_t(u)^{-1}) < 0
\end{align}
then by \ref{sec:local_linear_latent}, $u(t)$ is contracting to the subspace $\im(P)$ after a finite overshoot.

\subsection{Contraction to a submanifold}
\label{sec:manifold_contraction}
We now generalize contraction to submanifolds in $\R^n$ as explored in control contraction metrics (Section 5, \cite{manchester2017control}) to sufficient conditions for contraction to Hilbert and Banach submanifolds. 

We give a first result for Hilbert submanifolds, and then extend to Banach spaces using weighted SIPs as in \ref{sec:local_linear_latent}. As noted earlier, calculating the (weighted) contraction rate is often easier in $L^2$; thus, we can also combine the following results with $L^p$ comparison (\ref{sec:lp_estimates}) and a suitable embedding of a Hilbert submanifold $\mathcal{M}$ to obtain asymptotic convergence to submanifolds in $L^p$. We subsequently give the general result for Banach submanifolds.

\begin{theorem}[Contraction to submanifolds via weighted SIPs]
\label{thm:hbm_contraction}
Let $\phi \in C^\infty(V, W)$ be a total smooth map between Hilbert spaces $V, W$ such that $0$ is a regular value (see \cite{lang2012differential}) of $\phi$ (for all points in the zero set $u \in \phi^{-1}(0)$, $D\phi(u)$ is surjective). 
Let $\mathcal{M}  = \phi^{-1}(0)$, which is a time-invariant $\dim(W)$-submanifold of $V$. 

Define the weighted inner product:
\begin{align}
    \label{eq:manifold_weighted_ip}
    \langle x, y\rangle_{\phi} =: \langle D\phi(u)x, D\phi(u)y\rangle_W
\end{align}
as well as the associated Hilbert-space contraction rate $M^{\phi}$ \eqref{eq:l2_lognorm} (we drop the argument $u$ for conciseness, although \eqref{eq:manifold_weighted_ip} is smoothly $u$-varying).

Let $\dot{u} = f(t, u)$ be a differentiable vector field on $V$, and suppose that:
\begin{enumerate}
    \item $\forall u \in M, t \ge 0$, we have $D\phi(u)f(t, u) = 0$
    \item $u(t)$ is differentially contracting in the $\phi$-weighted inner product \eqref{eq:manifold_weighted_ip}:
    $$
    \sup_{t, u}M^\phi(Df_t(u)) = \lambda < 0
    $$
\end{enumerate}
Then $u(t)$ asymptotically approaches the manifold $\mathcal{M}$.
\end{theorem}
\begin{proof}
We show that $u(t)\to \mathcal{M} $ by showing that $\phi(u(t))$ is contracting to $0$. By hypotheses, we have:
\begin{enumerate}
    \item $\phi(u(t))$ has a fixed point at $0$, by either of the following equivalent conditions (compare to the linear version \eqref{eq:subspace_dynamics}):
    \begin{align}
    \label{eq:invariant_manifold}
    \forall u \in M, t \ge 0,\quad 
    \begin{cases}
        \frac{d}{dt}\phi(u) = D\phi(u)f(t, u) = 0\\
        f(t, u) \in T_u\mathcal{M}  = \ker(D\phi(u)) 
    \end{cases}
    \end{align}
    This states that the dynamics are tangent to $\mathcal{M}$, i.e. $\mathcal{M}$ is an invariant manifold of $v(t)$.

    \item Next, define $v = \phi(u)$; a perturbation $\delta u \in T_u V$ propagates as $\delta v = D\phi(u)\delta u$. Then, 
    \begin{align}
        \label{eq:manifold_projection}
        \frac{d}{dt}\frac12\norm{\delta v}_W^2 &= \Re \left\langle D\phi(u)\delta u, D\phi(u)Df_t(u)\delta u\right\rangle_W
    \end{align}
    
     By substitution of \eqref{eq:manifold_weighted_ip} into \eqref{eq:manifold_projection}, we have:
    \begin{align}
        \label{eq:manifold_contraction}
        \frac{d}{dt}\norm{\delta v}^2 &= \frac{2\Re \langle\delta u, Df_t(u)\delta u\rangle_{\phi}}{\langle \delta u, \delta u\rangle_{\phi}}\norm{\delta v}^2 \le 2M^{\phi}(Df_t(u)) \norm{\delta v}^2\\
        &\le 2\lambda \norm{\delta v}^2
    \end{align}
\end{enumerate}

Together, \eqref{eq:invariant_manifold} and \eqref{eq:manifold_contraction} imply that $\phi(u(t)) \to 0$ exponentially and thus $u(t)$ is contracting to the submanifold $\mathcal{M} = \phi^{-1}(0)$.
\end{proof}

\subsubsection{Banach submanifolds}
Result \eqref{eq:manifold_contraction} can also be extended to Banach spaces $V, W$ using weighted semi-inner products \eqref{eq:differential_weight} and the Dini identity \eqref{eq:dini_bound} analogously to section \ref{sec:local_linear_latent}. 

\begin{theorem}[Contraction to Banach submanifolds]
Let $\phi$ be as in \ref{thm:hbm_contraction} and define the weighted semi-inner product:
\begin{align}
\label{eq:bsp_weighted_norm}
(x, y)_+^{\phi} &:= (D\phi(u)x, D\phi(u)y)_+
\end{align}
letting $M^{\phi}$ be the associated weighted contraction rate. Let $\dot{u} = f(t,u)$ and suppose:
\begin{enumerate}
    \item $\forall u \in M, t \ge 0$, we have $D\phi(u)f(t, u) = 0$
    \item $u(t)$ is differentially contracting in the $\phi$-weighted \textit{semi}-inner product \eqref{eq:bsp_weighted_norm}:
    $$
    \sup_{t, u}M^\phi(Df_t(u)) = \lambda < 0
    $$
\end{enumerate}
Then $u(t) \to \mathcal{M} = \phi^{-1}(0)$.
\end{theorem}

\begin{proof}
 Let $u(t)$ be a solution, $v(t) = \phi(u(t))$, and $\delta u, \delta v$ be perturbations as before. Applying the Dini identity \eqref{eq:dini_bound} to the perturbation dynamics $\delta v(t)$, we have:
\begin{equation}
\begin{split}
    \label{eq:banach_manifold_contraction}
    D_t^+ \norm{\delta v} &= \frac{\left(\delta v, \frac{d}{dt}\delta v\right)_+}{\norm{\delta v}^2} \norm{\delta v} = \frac{\left(D\phi(u) \delta u, D\phi(u)Df_t(u)\delta u\right)_+}{\norm{D\phi(u)\delta u}^2} \norm{\delta v} \\
    &= \frac{\left(\delta u, Df_t(u)\delta u\right)_+^{\phi}}{\norm{\delta u}_{\phi(u)}^2} \norm{\delta v} \le M^{\phi}(Df_t(u))\norm{\delta v}
\end{split}
\end{equation}
Then by hypothesis (2), since $M^{\phi}(Df_t(u) < 0$ uniformly, $v(t)$ is contracting, and by hypothesis (1), $\mathcal{M}$ is an invariant manifold. Thus, $\phi(u(t)) \to 0$  and therefore $u(t)$ approaches the submanifold $\mathcal{M} = \phi^{-1}(0)$.
\end{proof}

Condition \eqref{eq:banach_manifold_contraction} can also be extended using weighted norms \eqref{eq:differential_weight}. We can apply an \textit{outer} (time-invariant) weight $\Theta(v)$ to the $\delta v$ dynamics via the weighted semi-inner product
$$
(x, y)_+^{\Theta, \phi} = (\Theta(v)D\phi(u)x, \Theta(v)D\phi(u)y)_+
$$
Then by results in section \ref{sec:local_linear_latent}, if the associated contraction rate is negative,
\begin{align}
    \label{eq:weighted_manifold_contraction}
    \sup_{t, u}M^{\Theta, \phi}(Df_t(u)) &< 0
\end{align}
then $v(t)$ is contracting to $0$ and $u(t)$ is contracting to $\mathcal{M}$ after a finite overshoot.

We note that both in the cases of contraction to subspaces \eqref{eq:subspace_contraction} and submanifolds \eqref{eq:banach_manifold_contraction}, the sufficient condition is a negative contraction rate in a weighted semi-inner product, with the primary difference being a uniform or space-varying weight.




\subsubsection{Order reduction, discovery of attractors} 
Suppose we have an autonomous system $\dot{v} = f(v)$, and we would like to investigate the existence of particular attractors. The submanifold contraction condition \eqref{eq:banach_manifold_contraction} suggests a discovery procedure for for any submanifold attractors via functional optimization. Let $W$ be a chosen target space of dimension $n$; then, if the following (nonconvex, nonlinear) problem (with $M^{\Theta, \phi}$ defined as in \eqref{eq:weighted_manifold_contraction}):
\begin{align*}
    \text{minimize}\ & \sup_{v \in V}M^{\Theta, \phi}\left(Df(v)\right) + \sup_{v \in V}\norm{D\phi(v)f(v)}\\
    \text{over}\ & \phi \in C^\infty(V, W), \Theta \in C^1(W, GL(W))\\
    \text{subject to } & 
        \rank(D\phi(0)) = n
\end{align*}
is feasible with objective less than $0$, the $n$-dimensional submanifold $\mathcal{M} = \phi^{-1}(0)$ is an attractor for the system. We can also replace minimization over the total space $V$ with some closed convex subset $W \subset V$ to establish $W$ as a basin of attraction. For work on related convex formulations for $\Theta$ in the finite-dimensional setting, we refer the reader to control contraction metrics \cite{manchester2017control}.

\section{Inheritance of contraction properties}
\label{sec:combinations}
We now extend the development of \textit{combination properties} (\cite{lohmiller1998contraction}, 3.8) of contraction rates of composite dynamical systems in Banach spaces, using properties of the semi-inner product such as subadditivity and positive scalar homogeneity (see \cite{lumer1961semi}).

\subsection{Additive combinations}
\begin{theorem}
Let two dynamical systems in a common Banach space $V$ be given by
\begin{align*}
    \dot{u}_1 = f_1(t, u_1),\quad \dot{u}_2 = f_2(t, u_2)
\end{align*}
with $f_1, f_2 \in C^1(\R^+\times V,V)$. Suppose that the systems are both contracting in the same weighted \eqref{eq:differential_weight} semi-inner product:
\begin{align*}
    \min(M^\Theta(Df_{1,t}(u_1)), M^\Theta(Df_{2,t}(u_2))) &= \lambda
\end{align*}
Then any linear combination of the form:
\begin{align*}
\dot{u} &= g(t,u) := \alpha_1(t)f_1(t, u) + \alpha_2(t)f_2(t, u)\\
0 &\le \sup_t \max(\alpha_1(t), \alpha_2(t)) \\
0 &< \epsilon := \inf_t (\alpha_1(t) + \alpha_2(t)) 
\end{align*}
is contracting.
\end{theorem}
\begin{proof}
Noting that sub-additivity of the semi-inner product (5.2, \cite{soderlind2006logarithmic}) carries over to the weighted SIP \eqref{eq:differential_weight}, the contraction rate of the composite system satisfies:
\begin{align*}
    M^\Theta(Dg_t(u)) &\le \alpha_1(t)M^\Theta(Df_{1,t}(u)) + \alpha_2(t)M^\Theta(Df_{2,t}(u)) \le 2\epsilon\lambda 
\end{align*}
\end{proof}

\subsection{Feedback combinations}
\label{sec:feedback_combinations}

\subsubsection{Skew-adjoint coupling in Hilbert spaces}
We first consider the case of a feedback coupling in Hilbert spaces $H_1, H_2$ having a \textit{skew-adjoint} form, generalizing the $\R^n$ case from (3.8.2, \cite{lohmiller1998contraction}). We endow the product space $H = H_1 \times H_2$ with an inner product induced by the sum: 
\begin{align}
\label{eq:product_euclidean}
\langle (f_1, f_2), (g_1, g_2)\rangle_H &= \langle f_1, g_1\rangle_{H_1} + \langle f_2, g_2\rangle_{H_2}
\end{align}
also known as the direct sum $H_1 \oplus H_2$. 

\begin{theorem}
Let two dynamics $u_1(t), u_2(t)$ in $H_1, H_2$ have the skew-adjoint feedback coupling:
\begin{equation}
\begin{split}
    \label{eq:skew_adjoint}
    \dot{u}_1 &= f_1(t, u_1, u_2)\\
    \dot{u}_2 &= f_2(t, u_1, u_1)\\
    J_{12}(t, u_1, u_2) &= -J_{21}(t, u_1, u_2)^*
\end{split}
\end{equation}
where $J$ is the Jacobian of $(f_1, f_2)$. The feedback system is contracting in $H$ if and only if both $u_1, u_2$ are contracting in $H_1, H_2$. 
\end{theorem}
\begin{proof}
Let $u = (u_1, u_2)$ and $f = (f_1, f_2)$ represent the composite system with dynamics $\dot{u} = f(t, u)$ in $H$. Its contraction rate is
\begin{equation}
\begin{split}
    \label{eq:skew_contraction}
    M(Df_t(u)) &= \sup_{\norm{v}_H = 1} \Re \left\langle J(t, u)v, v\right\rangle_H \\
    &= \sup_{\norm{v}_H=1} \Re\left\langle\begin{bmatrix}J_{11}(t, u) & J_{12}(t,u) \\ -J^*_{12}(t, u) & J_{22}(t,u)\end{bmatrix}\begin{bmatrix}v_1 \\ v_2\end{bmatrix}, \begin{bmatrix}v_1\\v_2\end{bmatrix}\right\rangle_H \\
    &= \sup_{v_1 \in H_1, v_2 \in H_2}\frac{\Re\left[\langle J_{11}(t,u) v_1, v_1\rangle_{H_1} +  \langle J_{22}(t,u)v_2, v_2\rangle_{H_2}\right]}{\langle v_1,v_1\rangle_{H_1} + \langle v_2, v_2\rangle_{H_2}} \\
    &= \max(M(J_{11}(t,u)), M(J_{22}(t,u)))
\end{split}
\end{equation}
Thus, the composite system is contracting if and only if the individual systems are contracting. 
\end{proof}




\subsubsection{Zero-range coupling in Banach spaces}
We now consider a generalization of a skew-adjoint feedback coupling \eqref{eq:skew_adjoint} to Banach spaces via the semi-inner product numerical range. 
\begin{definition}[Zero-range operator]
Recall that the \textit{numerical range} (Definition 4, \cite{lumer1961semi}) of an operator $A \in \B(V)$ associated with a semi-inner product $(\cdot, \cdot)_+$ on a Banach space $V$ (see \eqref{eq:sip}) is defined as 
$$
W(A) = \{(v, Av)_+\ |\ \norm{v} = 1\}
$$
We call $A$ \textit{zero-range} if $\Re[W(A)] = \{0\}$.
\end{definition}
Note that if $V$ is a Hilbert space, $A$ is a zero-range operator if and only if it is skew-adjoint. Unlike the Hilbert-space feedback coupling, we do not have an equivalence due to subadditivity of the semi-inner product, but two slightly different assumptions in each direction as we describe below.

Suppose we have $n$ dynamics in Banach spaces $V_1,..., V_n$ let and $u = (u_1,...,u_n)$:
\begin{align*}
    \dot{u}_i &= f_i(t, u)
\end{align*}
where the Jacobian of this feedback system is denoted by a matrix of bounded linear operators $J_{ij}(t, u) = \frac{\partial f_i}{\partial u_j}(t, u)$. Let $V = \prod_i V_i$ be given any $\ell^p$ norm 
\begin{align}
    \label{eq:lp_product}
    \norm{v}_V = \left(\sum_i\norm{v_i}_{V_i}^p\right)^{1/p},\quad 1 \le p \le \infty
\end{align} 

\begin{theorem}[Feedback contraction from individual contraction]
Let $u_1...u_n$ and $J$ be as above, and suppose the individual dynamics (i.e. the diagonal terms) are contracting:
\begin{equation}
\begin{split}
    \label{eq:decoupled_contraction}
    \forall i,\quad \lambda_i := \sup_{t, u \in V}M(J_{ii}(t, u)) < 0
\end{split}
\end{equation}
Then the feedback system $u(t)$ is contracting in $V$.
\end{theorem}
\begin{proof}
A perturbation $\delta u \in V$ in the feedback system satisfies the growth bound:
\begin{align*}
    \norm{\delta u(t)}_V &\le \left(\sum_i\norm{v_i}e^{p\lambda_it}\right)^{1/p} \le \norm{\delta u(0)}_V e^{\max_i\lambda_i t}
\end{align*}
hence the feedback system is contracting.
\end{proof}
In the opposite direction, we have:
\begin{theorem}[Individual contraction from feedback contraction]
Let $u_1...u_n$ and $J$ be as above, and decompose the Jacobian as $J = D + F$, where $D$ is the diagonal (self-exciting) components and $F$ is the off-diagonal (feedback) components. Suppose that:
\begin{enumerate}
    \item $F$ is \textit{zero-range}, i.e. $\Re (v, F(t, u)v)_+ = 0$ for all $t \ge 0$ and $u, v \in V$
    \item $M_V(D(t, u)) + M_V(F(t, u)) < 0$
\end{enumerate}
Then each system $u_i(t)$ is contracting in $V_i$.
\end{theorem}
\begin{proof}
Together (1) and (2) imply that \textit{both} the diagonal and composite system are contracting in $V$, i.e. $M_V(D(t,u)) < 0$ and $M_V(J(t,u)) < 0$. 

In the $\ell^p$ norm on $V$, the $V$-operator norm of $D$ is $\norm{D}_{\B(V)} = \max_i(\norm{D_{ii}}_{\B(V_i)})$. Hence, 
\begin{equation}
\begin{split}
    \label{eq:converse_max}
    0 &> M_V(D(t,u)) =  \lim_{h\to 0^+} \frac{1}{h}\left(\max_i\left(\norm{I + hD_{ii}(t, u)}_{\mathcal{B}(V_i)}\right) - 1\right) \\
    &= \max_i M_{V_i}(D_{ii}(t, u))
\end{split}
\end{equation}
thus, the individual systems are contracting \eqref{eq:decoupled_contraction}. 
\end{proof}
However, we note that a converse statement (that $F$ is zero-range if the individual systems and feedback system are contracting) is not possible in general due to sub-additivity of the SIP.

\subsection{Feedforward combinations}
\begin{theorem}
Let two dynamics in Banach spaces $V_1, V_2$ be:
\begin{align*}
    \dot{u}_1 &= f_1(t, u_1)\\
    \dot{u}_2 &= f_2(t, u_1, u_2)
\end{align*}
Suppose that:
\begin{enumerate}
    \item $\dot{u}_1$ is contracting: $\sup_{t, u_1}M\left(\frac{\partial f_1}{\partial u_1}(t, u_1)\right) = \lambda_1 < 0$
    \item The self-excitation of $u_2$ is contracting: $\sup_{t, u_1, u_2}M\left(\frac{\partial f_2}{\partial u_2}(t, u_1, u_2)\right) = \lambda_2 < 0$
    \item The feedforward coupling is uniformly bounded: $$
    \sup_{t, u_1, u_2} \norm{\frac{\partial f_2}{\partial u_1}(t, u_1, u_2)}_{\B(V_1, V_2)} = B < \infty$$
\end{enumerate}
Then the system $u_2(t)$ is contracting in $V_2$.
\end{theorem}
\begin{proof}
Applying the Dini identity \eqref{eq:dini_bound} to a perturbation $\delta u_2 \in T_{u_2}V_2$, we obtain the growth bound:
\begin{align*}
D_t^+ \norm{\delta u_2(t)} &\le \frac{(\delta u_2, \frac{\partial f_2}{\partial u_1}(t, u_1, u_2)\delta u_1)_+ + (\delta u_2, \frac{\partial f_2}{\partial u_2}(t, u_1, u_2)\delta u_2)_+}{\norm{\delta u_2}^2}\norm{\delta u_2}\\
&\le B \norm{\delta u_1(t)} + \lambda_2 \norm{\delta u_2(t)} \le Be^{\lambda_1 t}\norm{\delta u_1(0)} + \lambda_2 \norm{\delta u_2}\\
\implies \norm{\delta u_2(t)} &\le \norm{\delta u_2(0)}e^{\lambda_2 t} + \frac{B\norm{\delta u_1(0)}}{\lambda_1 + \lambda_2}e^{\lambda_1 t}
\end{align*}
where we used subadditivity and Cauchy-Schwarz of the semi-inner product (Definition 1, \cite{lumer1961semi}).
Thus, $u_2$ is contracting after a finite overshoot.
\end{proof}
This extends a same result for finite-dimensional systems with the $\ell^2$ norm (2.3, \cite{manchester2014combination}).

\subsection{Linear continuum couplings}
We now consider additive couplings on a continuum.

Let $V$ be a Banach space and $K \subset \R^n$ a closed $n$-cell with side length $\ell$. Let $u(t, x) : \R^+ \times K \to V$ be a continuum system of trajectories in $X$, with dynamics defined by:
\begin{align}
\label{eq:contsystem}
\frac{\partial}{\partial t} u(t, x) &= f(t, x, u)
\end{align}
with $f$ continuously differentiable. 
Let $v(t) \in V$ be defined as the linear combination
\begin{align}
\label{eq:cont_comb_0}
v(t) = \int_K \varphi(x)u(t, x)dx
\end{align}
where $\varphi : K \to \R$ is a weight function and $\int$ is the Bochner integral. 
\begin{theorem}[Continuum coupling with positive weight]
Let $u(t, x), v(t)$ be as in \eqref{eq:contsystem} and \eqref{eq:cont_comb_0}, and suppose that:
\begin{enumerate}
    \item $\varphi \ge 0$ on $K$
    \item $\int_K \varphi > 0$
    \item $f$ has the \textit{point-wise} contraction rate:
\begin{align*}
    \sup_{t, x, u}M\left(f_{t,x}\right) &= \lambda 
\end{align*}
where $f_{t,x}(u) = f(t,x,u)$. 
\end{enumerate}
Then, $v(t)$ is contracting in $V$.
\end{theorem}
\begin{proof}
Since $K$ has finite measure, we can exchange $\partial_t$ and $\int$ by Hille's theorem (see 1.19, \cite{van2008stochastic}) to obtain the $v$ dynamics:
\begin{align}
    \label{eq:cont_comb}
    \dot{v} = \int_K \varphi(x)f(t, x, u)dx
\end{align}
Then by sub-additivity and homogeneity of $M$ for nonnegative scalars,
\begin{align*}
    \sup_tM\left(\int_K \varphi(x)f_{t,x}dx\right) &\le \sup_t\int_K \varphi(x)M(f_{t, x})dx \le \lambda \int_K\varphi(x)dx = c\lambda
\end{align*}
for some $c > 0$. Thus, point-wise contraction implies combined contraction under the weight $\varphi$.
\end{proof}


\section{Symmetries in contracting dynamics}

We now consider properties of integral curves of vector fields which are contracting in some weighted semi-inner product and also possess some symmetry. These notions were developed by \cite{russo2011symmetries} for $\R^n$, which we extend in several ways to dynamics in Banach spaces.

\subsection{Spatial symmetries}
\label{sec:spatial_linear}
\begin{theorem}
Let $\dot{u} = f(t, u)$ in a Banach space $V$, and $\Gamma \le GL(V)$ a subgroup of bijective bounded linear operators on $V$ such that $f$ is $\Gamma$-equivariant in space:
\begin{align}
\label{eq:spatial1}
\forall T \in \Gamma, u \in V, t \ge 0,\quad f(t, Tu) = Tf(t, u)
\end{align}
If $f$ is contracting with respect to some weighted semi-inner product $M^\Theta$, then $u(t)$ approaches a single $\Gamma$-invariant vector exponentially fast.
\end{theorem}
\begin{proof}
Define the set of solutions at time $t$ as $S(t) = \{\phi(0, t, u)\ |\ u \in V\}$, where $\phi$ is the solution operator generated by $f$ (see \ref{sec:evolution_family}). By the equivariance \eqref{eq:spatial1}, we have
$$
\dot{u} = f(t, u) \implies T\dot{u} = Tf(t, u) = f(t, Tu)
$$
implying that $u \in S(t) \implies Tu \in S(t)$, thus $S(t)$ is $\Gamma$-invariant for all $t$. 

By the contraction hypothesis, $S(t)$ is then contracting after a finite overshoot to a single solution $u_*(t)$ which is $\Gamma$-invariant. 
\end{proof}

\subsubsection{Nonlinear symmetries}
\label{sec:nonlinear_symmetries}
The previous result extends straightforwardly to non-linear symmetries.

Suppose $\Gamma \subset C^1(X)$ is more generally a group of diffeomorphisms, such that $f$ is \textit{differentially} equivariant:
\begin{align}
\label{eq:nonlinear_equivariance}
\forall h \in \Gamma, u \in V, t \ge 0,\quad f(t, h(u)) = Dh(u)f(t, u)
\end{align}
(We note that this equivariance condition \eqref{eq:nonlinear_equivariance} also appears in an un-published manuscript by \cite{unpub_boffi}.)
Then, for any $h \in \Gamma$ we have $\frac{d}{dt}h(u) = Dh(u)\dot{u} = f(t, h(u))$, thus $S(t)$ is $\Gamma$-invariant. Finally, if $M^\Theta(f_t) < 0$ for some weight $\Theta$, 
$S(t)$ approaches a single solution $u_*(t)$ which is a fixed point of each $h \in \Gamma$.

\subsubsection{Contraction to a $\Gamma$-invariant subspace} 
Suppose we only wish to show that $u(t)$ approaches \textit{some} $\Gamma$-invariant solution rather than a unique one; we apply contraction in a weighted semi-inner product to do so.
\begin{theorem}
Let $\Gamma \le GL(V)$ and $\dot{u} = f(t, u)$ be $\Gamma$-equivariant as in \ref{sec:spatial_linear}. Define the $\Gamma$-invariant subspace:
\begin{align*}
    V_\Gamma = \{u \in V\ |\ Tu = u\ \forall T \in \Gamma\}
\end{align*}
and suppose that $\dim(V_\Gamma) < \infty$. Then, there exists some weight $\Theta$ such that if $M^\Theta(f_t) < 0$, solutions $u(t)$ approach $V_\Gamma$ exponentially fast.
\end{theorem}
\begin{proof}
By the Hahn-Banach theorem, $\dim(V_\Gamma) < \infty)$ guarantees the existence of a bounded linear projection operator $P$ with $\im(P) = V_\Gamma$ (see \cite{randrianantoanina2001norm}).

By $\Gamma$-equivariance, we have that:
\begin{align*}
    \forall T \in \Gamma, u \in V, t \ge 0,\quad Tf(t, Pu) &= f(t, TPu) = F(t, Pu)
\end{align*}
implying that $QF(t, Pu) = 0$ for all $t,u$ (where $Q = I - P$). Thus $V_\Gamma$ is a flow-invariant subspace of $u(t)$. 

Finally, applying the results for contraction to an invariant subspace (see \ref{sec:subspace_contraction}) contraction in the $Q$-weighted semi-inner product $M^{Q}(f_t) < 0$ implies that $u(t)$ is contracting to $V_\Gamma$.
\end{proof}

We can also extend this to (outer) $\Theta$-weighted SIPs as in \eqref{eq:weighted_subspace_contraction} to obtain contraction to $V_\Gamma$ after a finite overshoot.

\begin{example}[Periodic heat equation]
Let $\partial_{t} u = \Delta u$ on the torus $\T^n = \R^n/\Z^n$, and $\Gamma$ be a group of bounded linear translation operators. We have $\Delta T = T\Delta$ for all $T \in \Gamma$, and $\mu(\Delta) < 0$ in $L^2_0(\T^n)$, the set of mass-zero square-integrable functions (see Proposition 6.1, \cite{soderlind2006logarithmic} for the contraction rate); thus $u$ tends to a $\Gamma$-invariant solution. Since $\Gamma$ was arbitrary, this is an alternate proof that harmonics on $\T^n$ are translation-invariant.  
\end{example}




\subsection{Spatio-temporal symmetries}
\begin{theorem}
\label{thm:spatiotemporal}
Let $\Gamma \le GL(X)$ be a cyclic subgroup of order $k$ generated by a bijective bounded linear operator $T$, and $\dot{u} = f(t, u)$ with the $\Gamma$-spatiotemporal symmetry
\begin{align}
\label{eq:spatiotemporal}
\forall T \in \Gamma, u \in V, t \ge 0,\quad f(t, Tu) = Tf(t+\Delta t, u)
\end{align}
for some $\Delta t > 0$.
If $f$ is contracting in some weighted semi-inner product,
\begin{align}
\label{eq:spatiotemporal_contraction}
\sup_{t}M^\Theta(f_t(x)) = \lambda < 0
\end{align}
then solutions $u(t)$ asymptotically approach a $k\Delta t$-periodic function.
\end{theorem}
\begin{proof}
By the symmetry \eqref{eq:spatiotemporal}, we have that $f(t, u) = f(t + k\Delta t, u)$ for all $t, u$. Thus, if $u(t)$ is a solution, then letting $s = t + k\Delta t$,
$$
\dot{u}(s) = T^k\dot{u}(s) = T^kf(s, u(s)) = f(t, T^ku(s)) = F(t, u(s))
$$
implying that $u(s) = u(t + k\Delta t)$ is also a solution, i.e. $S(t)$ is $k\Delta t$-translation invariant (letting $S(t)$ be the set of solutions at time $t$ as in \ref{sec:spatial_linear}).

By the contraction hypothesis, $S(t)$ asymptotically approaches a single solution. Thus, for any $m, n \in \N$ we have:
$$
\norm{u(t + mk\Delta t) - u(t + nk\Delta t)} \le e^{\lambda \min(m, n)k\Delta t}\norm{u(t + k\Delta t) - u(t)}
$$
implying that $\{u(t + nk\Delta t)\}_n$ is a Cauchy sequence. Since $V$ is a Banach space and thus complete, this implies $\lim_{n\to\infty}u(t + nk\Delta t) = p$ for some $p \in V$. As this is true for arbitrary $t$, we have that solutions converge to a periodic function with period $k\Delta t$.
\end{proof}

\subsubsection{Nonlinear symmetries} 
\label{sec:nonlinear_spatiotemporal}
The previous result (\ref{thm:spatiotemporal}) extends to nonlinear spatiotemporal symmetries as follows.

\begin{theorem}
Let $h$ be a diffeomorphism on $V$ which generates a cyclic group of maps of order $k$ by the fixed point condition $(h \circ ... \circ h)(u) =: h^k(u) = u$ for all $u \in V$. Suppose there exists $\Delta t > 0$ such that for all $n \le k$, the dynamics satisfy the \textit{differential} spatiotemporal symmetry:
\begin{align}
    \label{eq:nonlinear_spatiotemporal}
    \forall u \in V, t \ge 0,\quad f(t, h^n(u)) = Dh^n(u)f(t + n\Delta t, u)
\end{align}
If $f$ is contracting in some weighted semi-inner product $M^\Theta(f_t) < 0$, then solutions tend exponentially to a $k\Delta t$-periodic function after a finite overshoot.
\end{theorem}
\begin{proof}
For any solution $u(t)$, letting $s = t + k\Delta t$, we have 
\begin{align*}
    \frac{d}{ds}u(s) = \frac{d}{ds}h^k(u(s)) = Dh^k(u(s))f(s, u(s)) = f(t, h^k(u(s))) = f(t, u(s)) 
\end{align*}
implying that $u(s) = u(t + k\Delta t)$ is also a solution. Applying an identical argument as in \ref{thm:spatiotemporal}, the contraction hypothesis then implies convergence of $u(t) \to u(t + k\Delta t)$.
\end{proof}

\subsubsection{Contraction to limit cycles}
\label{sec:cycle_contraction}
We now combine the criteria of differential equivariance with respect to non-linear symmetries \eqref{eq:nonlinear_spatiotemporal} and partial contraction to a submanifold \eqref{eq:manifold_contraction} to informally give generic criteria for convergence to a limit cycle in a Banach space. We do so by relaxing the condition of global contraction \eqref{eq:spatiotemporal_contraction} to that of contraction to a \textit{loop}.

Let $V$ be a Banach space and $\phi\in C^{\infty}(V, \R)$ a smooth map  such that $0$ is a regular value and $\phi^{-1}(0)$ is an embedding of the unit circle $S^1$ in $X$. For example, 
$$
\phi(f) = (f(0)^2 + f(1)^2 - 1)^2 + \norm{f}^2
$$ 
on $L^2([0, 1])$ (technically, $\phi^{-1}(0) \isoto S^1$ except on a set of measure zero). Denote this embedding by $\hat{S}^1 = \phi^{-1}(0)$. 

Let $\dot{u} = f(t, u)$ in $V$ and $h : V \to V$ be any diffeomorphism, and take the coordinate transform $v = h(u)$ with conjugate dynamics:
\begin{align}
\label{eq:projected_cycle}
\frac{d}{dt}v = Dh(u)\dot{u} = Dh(h^{-1}(v))f(t, h^{-1}(v)) =: g(t, v)
\end{align}
Suppose $v(t)$ satisfies the conditions for contraction to the submanifold  \eqref{eq:banach_manifold_contraction}  $\hat{S}^1$,
\begin{equation}
\begin{split}
    \label{eq:loop_contraction}
    &1.\ \forall t \ge 0, z \in \hat{S}^1,\quad D\phi(z)g(t,z)  = 0\\
    &2.\ \sup_{t, v \in V} M^\phi(Dg_t(v)) < 0
\end{split}
\end{equation}
then solutions $v(t)$ exponentially approach the unit circle $\hat{S}^1$. 

Next, we impose a temporal symmetry on the vector field at the loop; suppose that $g$ is $T$-periodic on $\hat{S}^1$ for some $T > 0$:
\begin{align}
    \label{eq:cycle_symmetry}
    \forall t\in \R^+,\ v \in \hat{S}^1,\quad  g(t, v) &= g(t+T, v)
\end{align}
which may be obtained by showing spatiotemporal symmetry \eqref{eq:spatiotemporal} with respect to a cyclic group of rotation operators  on $\hat{S}^1$. 

Then for any solution $v(t)$ of \eqref{eq:projected_cycle} with initial condition $v(0) \in \hat{S}^1$, we have $v_n(t) = v(t + nT)$ and thus $v_n(0) = v(nT)$ is also a solution for each $n \in \N$. However, this does not preclude solutions which are arc-wise periodic on $\hat{S}^1$. 

Thus, to obtain a limit cycle attractor we finally impose a nonvanishing condition:
\begin{align}
    \label{eq:cycle_orientation}
    \forall t\in \R^+,\ v \in \hat{S}^1,\quad  g(t, v) \ne 0
\end{align}
which by continuity of $g$ implies that solutions of \eqref{eq:projected_cycle} cannot change orientation nor accumulate at any point on $\hat{S}^1$.

Conditions \eqref{eq:cycle_symmetry}, \eqref{eq:cycle_orientation} imply that solutions starting in $\hat{S}^1$ are $T$-periodic cycles on $\hat{S}^1$ (with possibly non-unit winding number). Combined with the loop contraction condition \eqref{eq:loop_contraction}, this implies that solutions $v(t)$ starting anywhere in $V$ approach the limit cycle $\hat{S}^1$. 

Thus, the original dynamics $u(t)$ approach a limit cycle on the (possibly non-circular) loop $h^{-1}(\hat{S}^1)$.

\subsubsection{Period synchronization in heterogeneous limit cycle systems}
We now apply the condition for contraction to limit cycles in \ref{sec:cycle_contraction} to analyze phase-locking phenomena in heterogeneous coupled-oscillator systems. 

Let $V$ be a Banach space, $u \in V^n$, and 
\begin{equation}
\begin{split}
    \label{eq:hetero_sys}
    \dot{u}(t) &= f(t,u)
\end{split}
\end{equation}
be some feedback system (see \ref{sec:feedback_combinations}). We say that \eqref{eq:hetero_sys} is a \textit{period-synchronized} system of limit cycles if there exist diffeomorphisms $h_i$ such that their conjugate dynamics $v_i = h_i(u_i)$ as in \eqref{eq:projected_cycle} are contracting to $\hat{S}^1$ and have the temporal symmetry:
\begin{align}
\label{eq:common_period}
\forall i,\ t\in \R^+,\ v \in \hat{S}^1,\quad g_i(t, v) = g_i(t + T, v)
\end{align}
for some \textit{common} $T>0$. Since the diffeomorphisms $h_i$ are arbitrary, we can compose arbitrary rotations $r_i$ and time-dilations $d_i$ of $\hat{S}^1$ to restate \eqref{eq:common_period} as the equivalent (simpler) condition:
\begin{align}
    \label{eq:sync_cycle}
    \forall i,j,\ t\in \R^+,\ v \in \hat{S}^1,\quad g_j(t, v) = g_i(t,v) = g_i(t + T, v) 
\end{align}
That is, all systems tend to the \textit{same} $T$-periodic limit cycle on $\hat{S}^1$, and we may assume it is constant-speed.

Condition \eqref{eq:sync_cycle} is a point-wise property describing a system of limit cycles which is \textit{always} period-synchronized; we now apply contraction to submanifolds \ref{sec:manifold_contraction} and subspaces \ref{sec:subspace_contraction} to derive sufficient conditions for \textit{eventual} period-synchronization in \eqref{eq:hetero_sys}. 

Suppose there exist:
\begin{enumerate}
    \item Diffeomorphisms $h_i$ such that \textit{at least one} of the individual dynamics $h_i(u_i)$ is contracting to the submanifold $\hat{S}^1$ \eqref{eq:loop_contraction} and is also nonvanishing on $\hat{S}^1$ \eqref{eq:cycle_orientation}.
    
    \item Angles $\theta_i$ and a corresponding a phase-shift subspace $W$:
    \begin{align*}
        W &= \{[R(\theta_1)x, ..., R(\theta_n)x]\ |\ x \in X\} \subset V^n\\
        \im(P) &= W,\ P^2 = P
    \end{align*}
    (where $P$ is a bounded linear projection to $W$ and $R(\theta_i)$ are planar rotation operators preserving the set $\hat{S}^1$) such that the \textit{composite} dynamics on the torus $(\hat{S}^1)^n \isoto \hat{T}^n$ are contracting to $W$ (see \ref{sec:subspace_contraction}). In other words, the dynamics eventually have constant helicity on $\hat{T}^n$:
    \begin{align*}
        \sup_{t \ge 0,\ v \in \hat{T}^n} & M^Q(Dg(t, v)) < 0
    \end{align*}
    where $Q = I - P$ and $M^Q$ is defined as in \ref{sec:subspace_contraction}, with the product space $V^n$ given some suitable norm (e.g. $\ell^\infty$ as in \ref{sec:feedback_combinations}).
\end{enumerate}

Then, up to diffeomorphisms, at least one system $h_i(u_i(t))$ has $\hat{S}^1$ as a limit cycle, and all others converge to time-translations of it. This implies that the original dynamics \eqref{eq:hetero_sys} have periodic attractors tending to a common period $T$ (which may in general be time-varying).  

\section{Applications}
\label{sec:applications}

\subsection{Well-posedness of partial differential equations}
\label{sec:pde}
In the following, we use existence plus contraction to deduce some properties of partial differential equations.
\subsubsection{Regularity for initial-value problems} 
Consider a PDE of the form 
\begin{align}
    \label{eq:pde1}
    \partial_t u &= f(t, D^{\alpha^1} u, ..., D^{\alpha^n} u)
\end{align} where $u$ is defined on a spatial domain $\Omega = \R^d$ or $\Omega = \R^d / \Z^d$ and $\{\alpha^i\}_{i=1}^n$ are length-$k$ multi-indices. We use the mixed-partial derivative notation:
$$
D^{\alpha} = \frac{\partial^{|\alpha|}}{\partial x_1^{\alpha_1}...\partial x_d^{\alpha_d}},\quad |\alpha| = \sum_i\alpha_i
$$
We formally consider \eqref{eq:pde1}  as an ODE $\partial_t u = F(t, u)$ on the space of compactly supported smooth functions $C_c^\infty(\Omega)$, and suppose the existence of some smooth solution $\phi_*(t) : \R^+ \to C_c^\infty(\Omega)$ for all time. Let $M_{k, p}$ be the contraction rate in $C_c^\infty(\Omega)$ with respect to the Sobolev norm $W^{k,p}(\Omega)$. For example, in the case of $p = 2$, 
\begin{align}
\label{eq:kp_rate}
M_{k, 2}(F_t) = \sup_{\phi_1\ne\phi_2 \in C_c^\infty(\Omega)}\frac{\sum_{i=0}^k\Re \langle \phi_1^{(j)} - \phi_2^{(j)}, F(t,\phi_1^{(j)}) - F(t,\phi_2^{(j)})\rangle}{\sum_{j=0}^k\langle \phi_1^{(j)} - \phi_2^{(j)}, \phi_1^{(j)} - \phi_2^{(j)}\rangle}
\end{align}
and for $1 < p < \infty$, the semi-inner product is given uniquely by the Gateaux derivative of the norm (see \cite{zhang2009reproducing}), 
\begin{align*}
    (u, v)^{k, p} &= \sum_{j=0}^{k-1}(u^{(j)}, v^{(j)})^p
\end{align*}
such that the $k, p$-contraction rate admits a similar form to \eqref{eq:kp_rate}.
If $F$ has a maximal $k, p$-expansion rate $\lambda \in \R$,
\begin{align}
\label{eq:kp_expansion}
\sup_{t \ge 0} M_{k, p}(F_t) \le \lambda
\end{align}
then by the existence hypothesis, any smooth initial condition $\phi \in C_c^\infty(\Omega)$ grows in the Sobolev norm at most exponentially with rate $\lambda$, giving the following continuous dependence on initial conditions:
\begin{align}
\label{eq:regularity}
\norm{\phi(t) - \phi_*(t)}_{k,p} \le \norm{\phi(0) - \phi_*(0)}_{k, p}e^{\lambda t}
\end{align}
Thus, smooth initial conditions cannot lose their regularity in finite time. 

If furthermore, we define an extension of \eqref{eq:pde1} with $D^\alpha$ taken in the weak sense, i.e. $D^\alpha u := v$ where 
\begin{align*}
    \int_\Omega D^\alpha u \phi dx &= (-1)^{|\alpha|}\int_\Omega v D^\alpha \phi dx\quad \forall \phi \in C_c^\infty(\Omega)
\end{align*}
such that $\overbar{F}_t$  has domain $D(\overbar{F}_t) = W^{k, p}(\Omega)$, then if $\overbar{F}_t$ satisfies the expansion rate \eqref{eq:kp_expansion} in the weak sense, solutions to the modified equation cannot lose their weak-differentiability in finite time. 

As a corollary, we can always ``improve'' regularity by adding a term $g$ such that, for some finite $\lambda \in \R$,
\begin{equation}
\begin{split}
    \label{eq:regularized}
\partial_t u &= f(t, D^{\alpha^1} u, ..., D^{\alpha^n} u) + g(D^{\alpha^1} u, ..., D^{\alpha^n} u)\\
\lambda & \ge \frac{(u-v, F_t(u) + G(u) - F_t(v) - G(v))_+}{\norm{u-v}^2} 
\end{split}
\end{equation}
This can also be further extended via the use of weighted $k, p$-semi-inner products $M_{k,p}^\Theta$ as in \eqref{eq:differential_weight}, in which case the growth bound \eqref{eq:regularity} applies up to a multiplicative constant.

\textbf{Note on contractive-weak solutions.} The ``regularized'' equation \eqref{eq:regularized} suggests a notion of weak solution in the sense of a vanishing dissipative term, in the case where the worst-case contraction rate \eqref{eq:kp_rate} may be infinite. Suppose there exists some $g, G$ as in \eqref{eq:regularized} such that:
\begin{align*}
   \sup_{t \ge 0}\frac{(u-v, (F_t(u)- F_t(v)) + \epsilon (G(u)  - G(v)))_+}{\norm{u-v}^2} &\le \lambda(\epsilon)
\end{align*}
We take the limit as $\epsilon \to 0^+$, preserving the dissipating effect of $G$, and consider a ``contractive-weak'' solution to be the limit
\begin{align}
\label{eq:cont_weak}
u(t) \in C^0(\R^+ \times \Omega)\ |\ \lim_{\epsilon\to 0^+} u_\epsilon = u
\end{align}
where $u_\epsilon$ is a classical solution, whose existence is assumed, to the $\epsilon g(u)$-regularized equation \eqref{eq:regularized}. We consider the existence of this limit in some suitable norm, e.g. $L^\infty(\R^+\times \Omega)$. 

This notion of weak solution, which we informally mention and whose development we leave to later work, is one possible generalization of the vanishing viscosity method inspiring viscosity solutions introduced by \cite{crandall1992user}. In the  vanishing viscosity approach, $g$ is the Laplace operator $\Delta$, which is exponentially stable for fixed boundary conditions and dissipative for Dirichlet problems in $L^2$ (see Proposition 6.1, \cite{soderlind2006logarithmic}), while the general \textit{viscosity solution} $u_* \in C^0$ to degenerate-elliptic second-order equations of the form $F(x, u, Du, D^2u) = 0$ is characterized by a \textit{comparison} property for test functions of the form:
\begin{align}
\label{eq:visc_test}
\Phi^\pm(x_0) = \{\phi \in C^2(E)\ |\ \exists \epsilon > 0, \phi(x_0) = u_*(x_0), \pm\phi \ge \pm u_* \text{ on } N_\epsilon(x_0)\}
\end{align}
such that the vector field produces a monotonic (order-preserving) response:
\begin{align}
\label{eq:visc_comparison}
\pm F(x_0, \Phi^\pm(x_0), D\Phi^\pm(x_0), D^2\Phi^\pm(x_0)) \le 0
\end{align}
for all $x_0 \in E$ (Section 2, \cite{crandall1992user}). Similarly, while  \eqref{eq:cont_weak} provides a procedure for constructing weak solutions of equations possessing $n$-th order space derivatives, a comparison principle characterizing contractive-weak solutions without constructing an explicit sequence of classical solutions may exist.

The relationship between viscosity, contractive-weak, and weak solutions in the sense of distributions remains to be characterized, while the potential benefit of a contractive-weak formulation is its applicability to general nonlinear $n$-th order equations (via $n, p$-contraction rates). This proposal of contractive-weak solutions is related to prior work by \cite{evans1980solving} which constructed weak solutions by a sequence of accretive vector fields using an $L^\infty$ weak pairing. A further extension involves relaxing \eqref{eq:regularized} to contraction in some time/space-varying weighted semi-inner product space $M^\Theta$, giving contraction in the Sobolev norm after a finite overshoot.

\subsubsection{Existence and uniqueness for time-independent equations}
We now give sufficient conditions commonly used fixed-point arguments for showing existence and uniqueness of time-independent equations via Sobolev contraction rates.
Suppose we have a time-independent equation on $W^{k, p}(\Omega)$:
\begin{align}
\label{eq:time_indep}
f(u, D^{\alpha^1} u, ..., D^{\alpha^n} u) = 0
\end{align}
with derivatives $D^\alpha$ taken in the weak sense, and $\Omega = \R^n$ or $\R^n/\Z^n$. Suppose the map $F$ given by $u(x) \mapsto f(u, D^{\alpha^1} u, ..., D^{\alpha^n} u)(x)$ is well-defined and has domain $D(F) = W^{k,p}(\Omega)$, and has a negative $k,p$-contraction rate:
$$
M_{k,p}(F) = \sup_{u\ne v\in W^{k,p}(\Omega)}\frac{\sum_{j=0}^k(D^ju - D^jv, F(D^ju) - F(D^jv))_+}{\norm{u-v}_{k,p}^2} < 0
$$
Then the corresponding time-dependent autonomous equation induces a contraction mapping,
\begin{align*}
\partial_t u &= f(u, D^{\alpha^1} u, ..., D^{\alpha^n} u)\\
u &\mapsto \Phi(t, u)\\
u_* &\mapsto u_*
\end{align*}
and thus by the Banach fixed-point theorem, approaches a unique fixed point $u_*$ exponentially, which is also the unique solution to \eqref{eq:time_indep}. 

\begin{example}[Nonlinear Poisson equation]
Let $u : \Omega \to \R^n$ with $\Omega = [0, 1]^d$, and consider the second-order equation:
\begin{align}
    \label{eq:poisson_like}
    \Delta u + f(u) &= 0\\
    u &= 0 \text{ on } \partial \Omega
\end{align}
with $f \in C^1(\R^n, \R^n)$ and $\Delta$ defined as the continuous linear extension of the classical Laplacian densely $\nabla \cdot \nabla$ densely defined on $C_c^\infty(\Omega) \subset W_0^{1,2}$ (the Sobolev space $W^{1,2}$ with vanishing Dirichlet conditions). Applying the divergence theorem, we have:
\begin{equation}
\begin{split}
    \label{eq:laplace_neumann}
    M(\Delta) &= \sup_{\phi \ne 0}\frac{\langle \phi, \nabla \cdot \nabla \phi\rangle}{\langle \phi, \phi\rangle} = \sup_{\phi \ne 0}\frac{\int_{\Omega}(\nabla \cdot (\phi\nabla \phi) - \nabla \phi \cdot \nabla \phi)d\Omega}{\langle \phi, \phi\rangle} \\
    &= \sup_{\phi \ne 0}\frac{\int_{\partial \Omega}\phi\nabla \phi \cdot d\hat{\partial \Omega} - \langle \nabla \phi, \nabla \phi\rangle }{\langle \phi, \phi\rangle} = \sup_{\phi \ne 0}\frac{-\langle\nabla \phi, \nabla \phi\rangle}{\langle \phi, \phi\rangle} \le -\lambda(\Omega)
\end{split}
\end{equation}
for some $\lambda(\Omega) > 0$, by the Poincar\'e inequality. Thus, \eqref{eq:poisson_like} has unique solutions if
$$
M_2^\Theta(f) = \sup_{u\ne v \in \Omega}\frac{\Re\langle u-v,f(u)-f(v)\rangle}{\norm{u-v}^2} < \lambda(\Omega)
$$
by subadditivity of $M$, in \textit{any} invertible weight $\Theta$ \eqref{eq:differential_weight}. The typical uniqueness proof for the Poisson equation $\Delta u + f = 0$ does not apply to \eqref{eq:poisson_like}, since $f(u_1) \ne f(u_2)$ in general.
\end{example}

\subsection{Reaction-diffusion systems}
\label{sec:rd}
Let $\Omega = [0, 1]^d$, $V = W^{2,2}(\Omega)$, and consider a nonautonomous system of networked reaction-diffusion equations on $V$ with zero-flux boundary conditions:
\begin{equation}
\begin{split}
    \label{eq:rd}
    \partial_t u_i &= \alpha_i\Delta u_i + f_i(t, u_1, ..., u_n)\\
\hat{\partial \Omega} \cdot \nabla_x u_i & = 0
\end{split}
\end{equation}
where $\alpha_i > 0$ for all $i$ and $\hat{\partial \Omega}$ is normal to the boundary. (We take the operator $\Delta$ to be the continuous linear extension of the Laplacian $\nabla \cdot \nabla$ densely defined on $C_c^\infty(\Omega)$.) 

\textbf{Suppression of patterns.} As noted in \cite{aminzare2013logarithmic, aminzare2014guaranteeing}, the existence of a constant solution $u \equiv c$ to \eqref{eq:rd} along with a contractive reactive term $f$ implies that a Turing instability cannot exist, since the Laplacian with zero-flux conditions \eqref{eq:rd} is semi-contractive. However, we can further relax the global contraction condition $M(f) < 0$ to the condition of contraction to the subspace of constant functions $\{u \in L^2(\Omega)\ |\ u \equiv c \in \R\}$ by taking the $L^2$-orthogonal projection $P$ as in \ref{sec:subspace_contraction}. Letting $Q = I-P$, if
\begin{align*}
    Qf(t, Pv) &= 0\\
    \sup_{t \ge 0}M^{Q}(f_t) &< M^{Q}(\Delta)
\end{align*}
(with $M^Q$ defined as in \eqref{eq:subspace_contraction}) then pattern formation is suppressed. The first condition can be interpreted as the property that once synchronized, the system cannot de-synchronize.

\textbf{Excitement of patterns.} We now consider the opposite case, a continuum generalization of anti-synchrony explored by \cite{wang2005partial}, where the following two conditions are met by a two-component reaction-diffusion system. Let $P$ be the $L^2$-orthogonal projection onto $\Im(\Delta)$, and suppose:
\begin{enumerate}
    \item There exists $u_* \in \im(P), u_*\ne 0$ satisfying \begin{align*}
        f_1(t, u_*, -u_*) = - \alpha_1\Delta u_*\\
        f_2(t, u_*, -u_*) = -\alpha_2 \Delta u_*
    \end{align*}
    \item The additive combination is contracting in the $Q = I-P$-weighed semi-inner product:
    $$
    \sup_{t \in \R^+, u_1, u_2\in V}\left[M^{Q}(Df_{1,t}(u_1, u_2)) + M^Q(Df_{2,t}(u_1, u_2))\right] < (\alpha_1 + \alpha_2)\left|M^Q(\Delta)\right|
    $$
\end{enumerate}
Then, the combined dynamics $v = u_1 + u_2$ has a fixed point $v_* = 0$ by (1) and is contracting by (2) to the orthogonal complement of constant solutions. This implies a spatially inhomogeneous anti-synchronizing pattern with $u_1 \to -u_2$ exponentially in $\norm{\cdot}_V$.

\subsection{Scalar conservation laws}
Let $x \in \Omega$ with $\dim(\Omega) =d$ and $u \in W^{1, p}(\Omega, \C)$ and consider the partial differential equation in divergence form 
\begin{align}
\label{eq:div_form}
\partial_t u + \nabla \cdot f(u) &= 0
\end{align}
where $f$ is a vector field (also called the \textit{flux}) in $C^2(\R^+ \times \R, \R^d)$, and $\nabla = \nabla_x$ is the gradient with respect to the space variable $x$. This equation states that the first integral $\int_\Omega u$ is constant up to flux at the boundaries. The linearized dynamics of \eqref{eq:div_form} are given by a family of linear maps $A(u)$:
\begin{equation}
\begin{split}
    \label{eq:div_operator}
    \delta u &= A(u)\delta u \\
    A(u)v &= -\nabla \cdot f'(u)v = -(f''(u) \cdot \nabla u + f'(u) \cdot \nabla)v 
\end{split}
\end{equation}
In the previous section, we considered $f'(u) = -\nabla$, of which \eqref{eq:div_form} is a kind of generalization.
We first consider the case with $p = 2$ and $\Omega=\T^d$. Applying the divergence theorem,
\begin{align*}
    \langle v, A(u)v\rangle &= \langle v, (f''(u) \cdot \nabla u)v + f'(u) \cdot \nabla v\rangle = \int_\Omega (\nabla \cdot [v^2f'(u)] - v\nabla v \cdot f'(u)) d\Omega \\
    &= \int_{\partial \Omega} v^2f'(u) \cdot d\vec{\partial \Omega} - \langle v, f'(u) \cdot \nabla v\rangle = - \langle v, f'(u) \cdot \nabla v\rangle
\end{align*}
Consider a fixed-mass set $S(m) = \{u \in L^2(\Omega)\ |\ \int_\Omega u = m\}$. Each $S(m)$ is an invariant set of the conservation law \eqref{eq:div_form}, and $S(m') = S(m) + (m'-m)$, thus without loss of generality we consider the set of mass-zero functions $S(0)$. We consider contraction within these invariant sets by substituting the above derivation:
\begin{equation}
\begin{split}
    \label{eq:div_contraction}
    \sup_{u \in S(0)} M_2\left(A(u)\right) & = \sup_{u, v \in S(0)} \frac{\langle v, f'(u)\cdot \nabla v\rangle}{\langle v, v\rangle} 
\end{split}
\end{equation}
Generalizing the diffusion case \eqref{eq:laplace_neumann}, we can conclude contraction in \eqref{eq:div_contraction} for $f'(u) = (\nabla)^k$ where $k$ is odd and vanishing boundary conditions of all orders, using integration by parts and the fact that $0$ is the only constant function in $S(0)$. 

We note that while \eqref{eq:div_contraction} gives a sufficient condition for contraction in $L^2$, often a more natural contraction metric for transport equations is the quadratic Wasserstein distance $W_2(u, v) = \left(\inf_{\pi \in \Gamma(u, v)}\int_{\Omega^2}\norm{x-y}^2d\pi(x,y)\right)^{1/2}$. For the linearization \eqref{eq:div_operator}, it is known (see \cite{peyre2018comparison}, equation (5)) that $W_2(u, u+\delta u)$ is well-estimated in the special case of positive densities $u$ by the $u$-weighted homogeneous Sobolev norm $\norm{\delta u}_{\dot{H}^{-1}(u)}$, in which case one can develop contraction in $W_2$ by a compatible semi-inner product \eqref{eq:sip}. Indeed, Wasserstein-norms generalizing $W_1$ to signed measures of hetereogeneous masses and contraction rates in these norms for nonconservative transport and reaction-diffusion systems are an area of active research (see \cite{piccoli2019wasserstein}).

\subsection{Functional regression in a Banach space}
\label{sec:regression}
Recent work in kernel methods for machine learning by \cite{zhang2009reproducing, zhang2012regularized, lin2019reproducing} have extended representer-type theorems from Hilbert spaces to Banach spaces via the Gateaux semi-inner product \eqref{eq:sip}. Zhang \textit{et al} show that $L^p(\Omega, \C)$ for $1 < p < \infty$ (more generally, any Gateaux-differentiable, uniformly convex, and dual-uniformly convex, per Theorem 9 of \cite{zhang2009reproducing}) admits a unique \textit{reproducing kernel} $K : \Omega^2 \to \C$ such that:
\begin{align}
    \label{eq:repro}
    \forall x \in \Omega, u \in V,\quad u(x) &= (u, K(x, \cdot))_V
\end{align}
where the RHS of \eqref{eq:repro} is also referred to as the \textit{evaluation functional} $E_x(u)$ in $V$, and we use $(u, v)_V := (u, v)_+ = (u, v)_-$ to denote the Gateaux semi-inner product on $V$. 
Let $\ell(x, y) : \Omega^2 \to \R^+$ be a convex loss function, $\{(x_i, y_i)\}$ a set of samples of the target function and $L(u)$ the empirical risk functional
\begin{align}
    \label{eq:fr}
    L(u) &= \sum_{(x_i, y_i)}\ell(u(x_i), y_i)
\end{align}
Then \cite{zhang2012regularized} note that by assumption of primal- and dual-uniform convexity, $L$ is Fr\'echet-differentiable and the functional gradient is given by 
\begin{align}
    \label{eq:grad_fr}
    DL(u) &= \sum_{(x_i, y_i)} \frac{\partial \ell}{\partial x}(u(x_i), y_i)K(x_i, \cdot)^*
\end{align}
where $K(x_i, \cdot)^*$ is given by the duality pairing $K(x_i, \cdot)^*(v) = (v, K(x_i, \cdot))_V$ (see Theorem 6, \cite{giles1967classes} for the existence and uniqueness of dual pairings arising from uniform convexity of $V$). Theorem 8 in the same study shows that $DL(u_*) = 0$ in \eqref{eq:grad_fr} if and only if $u_*$ is the minimizer of \eqref{eq:fr}. We define the resulting functional regression using the \textit{mirror descent}:
\begin{equation}
\begin{split}
    \label{eq:mirror_descent}
    u^*(v) &= (v, u)_V\ \forall v \in V\\
    \frac{d}{dt}u^* &= -\alpha DL(u)
\end{split}
\end{equation}
As the authors note in the conclusion, the convexity of the characterization equation $DL(u) = 0$ (and stability of the functional regression) presents a new challenge over the RKHS version of \eqref{eq:grad_fr} due to the nonlinearity and subadditivity of the semi-inner product.

Here, we present a simple sufficient condition for asymptotic convergence of $u(t) \to u_*$ using contraction in arbitrarily weighted SIPs \eqref{eq:differential_weight}, similar to contraction analysis of natural gradient descent in $\R^n$ presented by \cite{wensing2020beyond}. Let $\Theta(t, u) \in GL(V^*,V^*)$ be any uniformly bounded family of invertible operators as in \ref{sec:local_linear_latent},  $M^\Theta$ be the contraction rate with respect to the weighted SIP $(\cdot, \cdot)_V^\Theta$, and define 
$$
H(u) := \frac{\partial}{\partial u^*}DL(u)
$$
where in the case $V = \ell^2(\R^n)$, $H(u)$ is the Hessian. 
Then, if \begin{align*}
    M^\Theta(-H(u)) &= \sup_{v^* \in V^*}\frac{\Re (\Theta(t, u)v^*, (\dot{\Theta_t}(u) - D\Theta_t(u)H(u))\Theta(t,u)^{-1}v^*)_{V^*}}{\norm{\Theta(t,u)v^*}_{V^*}^2}\\
    &= \lambda < 0
\end{align*}
since $u_*$ is an equilibrium of \eqref{eq:mirror_descent} by Theorem 8,  \cite{zhang2012regularized}, we have that $u(t) \to u_*$ in $V$ exponentially with rate $\lambda$ after a finite overshoot. 


\section{Conclusion}
In this work, we developed several generalizations of foundational results in contraction theory, from weighted norms on the tangent space (\cite{lohmiller1998contraction}), to contraction to invariant subspaces (\cite{pham2007stable}) and submanifolds, and contraction with symmetric vector fields (\cite{russo2011symmetries}) to normed spaces lacking an inner-product structure. These results followed essentially by extension of the classical development of semi-inner product spaces (\cite{lumer1961semi}) to weighted spaces which are induced by objects of interest such as functions whose zero-sets are manifolds or bounded linear projections to subspaces. We showed how asymptotic properties implied by contraction in weighted spaces are related to other dynamical invariants such as Lyapunov exponents, yet uniform contraction rates are heavily dependent upon the norm. Using contraction in weighted semi-inner products, we derived conditions for more general asymptotic properties for dynamical systems in normed spaces, such as convergence to submanifolds limit cycles. Lastly, we introduced applications of this theory to the analysis of PDEs, including regularity, uniqueness, and stability of parabolic and transport equations.

\textbf{Future work.} Many physical evolution equations such as those found in quantum mechanics, diffusion/transport, and delay-differential equations admit descriptions as dynamical systems on Banach spaces; the choice of particular weighted semi-inner product spaces for establishing asymptotic stability of such systems is an area of future work. In particular, development of contraction in Wasserstein-like norms by \cite{piccoli2019wasserstein} for signed measures seems like a promising avenue for stability analysis of systems of transport equations. Additionally, adapting contraction theory for dynamical systems in Banach spaces to PDEs requires a systematic approach for extending the domain of linear and nonlinear differential operators in some reasonable way to a complete vector space, as discussed in \ref{sec:ode_pde}. Such methods would enable broad applicability of contraction theory to control and observer design for PDEs, such as decoherence in Lindblad dynamics describing coupled quantum systems (see \cite{rouchon2013contraction}) and nonlinear transport equations such as the Navier-Stokes (see \cite{coron2019controllability}). Furthermore, the relationship between contraction in smoothly-weighted semi-inner product spaces and contraction of dynamics on true infinite-dimensional manifolds (such as statistical manifolds, loop spaces, and infinite-dimensional Lie groups) remains an open question. This relationship is significantly complicated by the fact that, unlike the Hilbert/Riemannian case, multiple notions of geodesic distance are possible, as discussed in \cite{balestro2017angles}. While classical results in Hilbert manifold theory by \cite{henderson1969infinite} suggest that the existence of homeomorphisms to open sets of the model space in the separable infinite-dimensional setting may enable one to establish asymptotic convergence of solutions by showing contraction in some metric on $L^2$, a clear theory for dynamics on Banach or Fr\'echet manifolds, in the vein of work by \cite{simpson2014contraction} for finite-dimensional systems remains to be developed. Lastly, contractive dynamics as a regularity-controlling mechanism for PDEs suggests several types of weak solutions and methods for growth estimation.
    
\section*{Acknowledgements}    
This research did not receive any specific grant from funding agencies in the public, commercial, or not-for-profit sectors.

\appendix

\section{Nonlinear semigroups and evolution families}
\label{sec:pointfree}
In this section, we discuss the relationship between the upper Dini derivative of the norm of a perturbation and its relationship to the \textit{logarithmic norm} of the linearized vector field, show how it arises in much the same way as the Dini identity for semi-inner products \eqref{eq:dini_bound}. Theorem \ref{thm:orig_contraction} illustrates the relationship between ``pointfree'' methods from functional analysis such as the spectral properties of one-parameter semigroups, and the ``point-wise'' analysis used in the preceding sections.

Many vector-valued dynamical systems arising from well-posed initial value problems, such as ODEs and PDEs, admit natural representations as families of nonlinear operators on Banach spaces (see \cite{crandall1972nonlinear}). We briefly introduce two such algebraic structures corresponding to autonomous and nonautonomous dynamics in Banach spaces.
\begin{definition}[One-parameter semigroups of autonomous flows]
Consider the initial-value problem
\begin{align}
    \label{eq:cauchy_autonomous}
    \dot{u} &= f(u),\quad u(0) = u_0
\end{align}
whose existence of differentiable solutions $u(t)$ on some interval $[0, T]$ and uniqueness for any initial condition $u_0 \in V$ we assume. We represent the flows generated by $f$ by a \textit{one-parameter semigroup} of solution operators, i.e. a map $\Phi(t, u) : \R^+ \times V \to V$ satisfying the following properties:
\begin{enumerate}
    \item (Propagator) $\Phi(0, u_0) = u_0, \quad \Phi(t, u_0) = u(t)$
    \item (Semigroup) For all $s, t \ge 0, u \in V$, we have $\Phi(t, \Phi(s, u)) = \Phi(t + s, u)$. 
    \item (Continuity) For all $u \in V$ and $t \in \R^+$, we have $\norm{\Phi(s, u) - \Phi(t, u)} \to 0$ as $s \to t^+$.
    \item (Infinitesimal generator) For all $u \in V$, we have $\lim_{t \to 0^+}\norm{\frac{\Phi(t, u) - u}{t}} = f(u)$.
\end{enumerate}
\end{definition}
For non-autonomous systems, the solution operators $\Phi$ are in general dependent on the starting time; thus, $\Phi$ has a \textit{cocycle} rather than a semigroup structure (see \cite{latushkin1999evolution}).
\begin{definition}[Cocycles for non-autonomous systems]
Consider a non-autonomous initial-value problem with the same existence and uniqueness properties as \eqref{eq:cauchy_autonomous}:
\begin{align}
    \label{eq:nonautonomous}
    \dot{u} &= f(t, u), \quad u(0) = u_0
\end{align}
with $f \in C^1(\R^+ \times V, V)$. Then $f$ generates a solution map $\Phi(s, t, u)$ with the properties:
\begin{enumerate}
    \item (Propagator) For all $s \le t \in \R^+, u \in V$, we have $\Phi(t, t, u) = u, \Phi(s, t, u(s)) = u(t)$ 
    \item (Cocycle) For all $0 \le r \le s \le t, u \in V$, we have $\Phi(r, t, u) = \Phi(s, t, \Phi(r, s, u))$
    \item (Continuity) For all $u \in V$ and $t_0 \le t \in \R^+$, we have $\norm{\Phi(t_0, s, u) - \Phi(t_0, t, u)} \to 0$ as $s \to t^+$ 
    \item (Infinitesimal generator) For all $u \in V$ and $s \in \R^+$, we have:
    $$
    \lim_{t \to s^+}\norm{\frac{\Phi(s, t, u) - u}{t-s}} = f(s, u)
    $$
\end{enumerate}
\end{definition}

\section{Logarithmic norm and stability of nonlinear evolution families}
The following theorem, synthesized from several results in Chapter 5 of \cite{ladas1972differential}, establishes the central differential inequality relating the logarithmic norm of the linearized vector field (infinitesimal generator) and the contraction rate of its generated cocycle. 
\begin{theorem}[Vector fields with negative logarithmic norm are contracting]
\label{thm:orig_contraction}
Let $V$ be a Banach space, $f(t, u) \in C^1(\R^+ \times V, V)$ a differentiable time-varying vector field, and 
$$
\dot{u} = f(t, u)
$$

If for all $t \in [0, \infty)$ and $u \in V$, the Fr\'echet derivative $Df_t(u) \in \mathcal{B}(V)$ satisfies a \textit{logarithmic norm bound}:
$$
\mu(Df_t(u)) := \lim_{h\to0^+}\frac{\norm{I + hDf_t(u)}_{\mathcal{B}(X)} - 1}{h} \le \lambda(t)
$$
for some $\lambda(t) \in C([0, \infty))$, then for any initial conditions $u, v \in V$ and $t \ge s$, the evolution family is contracting:
$$
\norm{\Phi(t, s, u) - \Phi(t, s, v)}_V \le e^{\int_{s}^t\lambda(t
')dt'}\norm{u - v}_V
$$
or equivalently, $\Phi$ has bounded Lipschitz constant $L[\Phi(t, s, \cdot)] \le e^{\int_{s}^t\lambda(u)du}$. If furthermore $\sup_t\lambda(t) = -\alpha < 0$ for some $\alpha > 0$, the system is contracting with rate $\alpha$.
\begin{proof}
Let $0 \le s < t$. We start by applying the fundamental theorem of calculus to Fr\'echet spaces, 
$$
\Phi(t, s, u) = \Phi(t, s, v) + \int_0^1D\Phi_{t, s}(u + \xi(v-u))(v-u)d\xi
$$
Thus we obtain an upper bound on the distance at time $t$,
$$
\norm{\Phi(t, s, u) - \Phi(t, s, v)}_V \le \norm{u - v}_V\sup_{\xi \in [0, 1]} \norm{D\Phi_{t, s}(u + \xi(v-u))}_{\mathcal{B}(V)}
$$
Let $w := x + \xi(v-u)$. By the chain rule, the inner term has dynamics, defined for $t > s$,
$$
\frac{\partial}{\partial t}\left[D\Phi_{t, s}(z)\right] = Df_t(z)D\Phi_{t, s}(z) 
$$
Let $U_z(t) := D\Phi_{t, s}(z) \in \mathcal{B}(X, X)$, so that $\dot{U}_z = Df_t(z)U_z$. Now we obtain an upper bound estimate for $\norm{U_z(t)}$ as follows. Let $h > 0$; then 
\begin{align*}
    \norm{U_z(t + h)} &\le \norm{U_z(t) + hDf_t(z)U_z(t)} + o(h) \implies \\
    \frac{\norm{U_z(t+h)} - \norm{U_z(t)}}{h} &\le \frac{\norm{U_z(t)}\norm{I + hDf_t(z)} - \norm{U_z(t)} + o(h)}{h}
\end{align*}
Taking the limit as $h \to 0^+$, we have:
$$
D_t^+\norm{U_z(t)} \le \norm{U_z(t)}\lim_{h\to0^+} \frac{\norm{I + hDf_t(z)} - 1}{h} = \norm{U_z(t)}\mu(Df_t(z))
$$
where $D_t^+$ indicates the upper Dini derivative. Then by hypothesis, 
$$
D_t^+\norm{U_z(t)} \le \lambda(t)\norm{U_z(t)}
$$
Applying a Gr\"onwall-type inequality for Dini derivatives (see Lemma 11, \cite{davydov2021non}), we have 
$$
\norm{U_z(t)} \le \norm{U_z(s)}e^{\int_{s}^t\lambda(u)du} 
$$
where $\norm{U_z(s)} = \norm{D\Phi_{s, s}(z)} = \norm{DI(z)} = 1$.
Substituting into the original estimate, this gives:
\begin{align*}
\norm{\Phi(t, s, u) - \Phi(t, s, v)}_V &\le \norm{u - v}_V\sup_{\xi \in [0, 1]}\norm{U_{u + \xi(v-u)}(t)}_{\mathcal{B}(V)}\\ 
&\le e^{\int_{s}^t\lambda(t')dt'}\norm{x - y}_V
\end{align*}
thus if $\sup_t \lambda(t) = -\alpha < 0$, $\Phi_{t, s}$ is a contracting evolution family with exponential rate $\alpha$. 
\end{proof}
\end{theorem}

\section{Trace conditions for contraction in $\C^n$}
We remark on a fact which holds for feedback systems in finite-dimensional complex vector spaces due to convexity of the numerical range of the Jacobian. Let $\dot{u} = f(t, u)$ in $\C^n$ and $J(t, u)$ be its Jacobian as before. 

Suppose that $\Tr J(t, u) = 0$, that is, the vector field $f$ is divergence-free:
\begin{align*}
    \nabla \cdot f(t, u) &= 0
\end{align*}
A corollary of the Toeplitz-Hausdorff theorem as described in \cite{shapiro2004notes} is that $J(t, u)$ is then unitarily equivalent to a matrix $F(t, u)$ whose diagonal is zero, that is, there exists unitary $U(t, u)$ such that 
\begin{align*}
    U(t, u)^*J(t, u)U(t, u) &= F(t, u)
\end{align*}
Thus if $f$ is divergence-free and the feedback coupling under this unitary equivalence is contracting with rate $M(U^*JU) < 0$, the feedback system is contracting.

As a further corollary, if $f$ has uniformly negative divergence:
$$
\nabla \cdot f(t, u) = -\frac{c}{n},\quad c > 0
$$
Then letting a trace-zero modified Jacobian be $\overbar{J} = J + cI$, we have:
\begin{align*}
    M(J) \le M(\overbar{J}) - c = M(U^*\overbar{F}U) - c 
\end{align*}
Hence the system is contracting if the feedback coupling under a unitary equivalence has maximal expansion rate $c$.

\section{Contraction to invariant sets}
We extend the condition for contraction to zero-sets which are submanifolds \eqref{eq:manifold_contraction} to slightly more general zero-sets of a $C^1$ function $\phi$. Let $E = \phi^{-1}(0)$ (by continuity,  $E$ must be closed). We assume that $E$ is invariant for dynamics $\dot{u} = f(t, u)$, that is
\begin{align*}
    \forall t \ge0, u \in E,\quad D\phi(u)f(t, u) = 0
\end{align*}
Note that if $E$ is not path-connected, by regularity of $\phi$ there always exists a critical point $u_* \not\in E$ such that $D\phi(u_*) = 0$, in which case $\sup_{t, u}M^{\phi(u)}\left(Df_t(u)\right) \ge 0$. These are unstable equilibria. However, if we suppose that this set is of measure zero, i.e.
$$
\forall t\ge 0,\quad m(f_t^{-1}(0) \setminus E) = 0
$$
and verify the essential supremum contraction rate:
$$
\esssup_{t, u}M^{\phi}(Df_t(u)) < 0
$$
Then the flow has no attractors outside $E$ (of positive measure) and initial conditions almost everywhere are contracting to $E$ in the set distance $d(u(t), E)$ induced by the norm.

\bibliographystyle{elsarticle-num-names} 
\bibliography{references}

\end{document}